\documentclass[journal]{IEEEtran}

\ifCLASSINFOpdf

\else

\fi

\usepackage[comma,numbers,square,sort&compress]{natbib}
\usepackage{epstopdf}
\usepackage[nooneline,flushleft]{caption2}
\usepackage{mathrsfs}
\usepackage{graphicx,subfigure}
\usepackage{color}
\usepackage{amsfonts}
\usepackage{amssymb}
\usepackage[all]{xy}
\usepackage{bm}
\usepackage{multirow}
\usepackage{algorithm}
\usepackage{algorithmic}
\usepackage{amsthm}
\usepackage{amsmath}
\usepackage{cases}
\usepackage{fancyhdr}

\newtheorem{corollary}{Corollary}

\newtheorem{assumption}{Assumption}
\newtheorem{remark}{Remark}
\newtheorem{lemma}{Lemma}
\newtheorem{theorem}{Theorem}
\newtheorem{example}{Example}

\newcommand{\BEASN}{\begin{eqnarray*}}
\newcommand{\EEASN}{\end{eqnarray*}}
\newcommand{\BEAS}{\begin{eqnarray}}
\newcommand{\EEAS}{\end{eqnarray}}
\newcommand{\BEQ}{\begin{equation}}
\newcommand{\EEQ}{\end{equation}}
\newcommand{\BIT}{\begin{itemize}}
\newcommand{\EIT}{\end{itemize}}

\newcommand{\ie}{{i.e.}}
\newcommand{\nn}{{\nonumber}}

\newcommand{\calO}{\mathcal{O}}

\newcommand{\vc}{\textcolor[rgb]{0.00,0.00,0.00}}

\begin{document}
%

\title{Push-sum Distributed Dual Averaging for Convex Optimization in Multi-agent Systems with Communication Delays}

%
%
%

\author{Cong Wang, Shengyuan Xu, Deming Yuan, Baoyong Zhang, Zhengqiang Zhang
\thanks{C. Wang, S. Xu, D. Yuan and B. Zhang are with the School of Automation, Nanjing University of Science and Technology, Nanjing 210094, China (e-mail: jsyzwangcong@163.com; syxu@njust.edu.cn; dmyuan1012@gmail.com; baoyongzhang@njust.edu.cn).}\\

\thanks{Z. Zhang is  with the School of Electrical Engineering and Automation, Qufu Normal University, Rizhao 276826, China (e-mail: qufuzzq@126.com).}
}

\maketitle

\thispagestyle{fancy}
\chead{\bfseries\large This work has been submitted to the IEEE for possible publication. Copyright may be transferred without notice, after which this version may no longer be accessible.}
\cfoot{\quad}
\renewcommand{\headrulewidth}{0.4pt}
\renewcommand{\footrulewidth}{0pt}

\begin{abstract}
The distributed convex optimization problem over the multi-agent system is considered in this paper, and it is assumed that each agent possesses its own cost function and communicates with its neighbours over a sequence of time-varying directed graphs. However, due to some reasons there exist communication delays while agents receive information from other agents, and we are going to seek the optimal value of the sum of agents' loss functions in this case. We desire to handle this problem with the push-sum distributed dual averaging (PS-DDA) algorithm which is introduced in \cite{Tsianos2012}. It is proved that this algorithm converges and the error decays at a rate $\calO\left(T^{-0.5}\right)$ with proper step size, where $T$ is iteration span. The main result presented in this paper also illustrates the convergence of the proposed algorithm is related to the maximum value of the communication delay on one edge. We finally apply the theoretical results to numerical simulations to show the PS-DDA algorithm's performance.

\end{abstract}
\begin{IEEEkeywords}
 Distributed optimization, multi-agent system, push-sum, distributed dual averaging, communication delays.
\end{IEEEkeywords}

\section{Introduction}
With the development of wired and wireless technology, large-scale networks have been paid a large amount of attention to, including machine learning \cite{Liu2018,Lee2018,Bianchi2016,Mateos-Nunez2017}, Internet \cite{Hosseini2016}, game theory \cite{JLi2018},  mobile \cite{Demetriou2010,Zhang2018} and many other fields. All these systems can be considered as multi-agent systems with thousands of agents. In a large-scale network, it is shared to have communication delays which may affect computation results and cause some calculation errors. Then, we would like to study the convex optimization problems in multi-agent systems with interaction delays.

Due to the technology change of industrial engineering, demands for large-scale computation have been inevitable and it is hard and inefficient to deal with such problems with centralized optimization algorithms. Then utilizing distributed optimization algorithms which can calculate with various nodes in a network to tackle the convex optimization problems is a natural development trend. A great deal of researchers have proposed many useful algorithms for distributed optimization problems \cite{nedic2015,Yuan2021,Liu2021,Li2018,Nedic2010,Doan2019,Yuan2018,Xiong2021,Duchi2012,Yuan2014,Nedic2016,Akbari2017,JLi2020}. For distributed optimization algorithms, gradient descent methods are very common \cite{Nedic2010,Yuan2021,Liu2021,Li2018,nedic2015,Nedic2016}. \cite{Nedic2010} proposes a distributed subgradient algorithm to dope out solutions to the distributed convex optimization problems having nondifferentiable convex objective functions. Based on some earlier researches, the mirror descent algorithms are developed and have applications to distributed convex optimization problems \cite{Doan2019,Yuan2018,Xiong2021}. In \cite{Doan2019}, the authors analyze centralized and distributed mirror descent algorithms, and \cite{Yuan2018} establishes the convergence results of the distributed stochastic mirror descent algorithm and the epoch-based distributed stochastic mirror descent algorithm for the constrained distributed optimization problem with strongly convex objective functions. Recently, dual averaging algorithms have a great number of applications \cite{Duchi2012,Yuan2014} for the splendid convergence performance and simplified convergence analysis. Authors of \cite{Duchi2012} investigate a dual averaging subgradient method, and authors of \cite{Yuan2014} discuss an inexact version of dual averaging algorithm.

Push-sum method is deemed to one efficient protocol to tackle consensus problems for multi-agent system. Various researchers have developed some new algorithms by combining the push-sum method and some other algorithms. \cite{nedic2015} and \cite{Nedic2016} introduce an algorithm called the subgradient-push method, which is constructed by the push-sum method and the gradient descent protocol. This algorithm is practical for convex optimization problems with time-varying directed communication graphs. In \cite{Tsianos2012}, authors formulate the PS-DDA algorithm that constructed by push-sum averaging method and distributed dual averaging protocol. \cite{Lee2018} proposes an online algorithm which consists of online dual averaging algorithm and push-sum protocol for solving the decentralized online convex optimization problems. In \cite{Yuan2015}, authors develop a push-sum gradient-free method using the push-sum algorithm and a gradient-free step. Motivated by the subgradient-push method in \cite{nedic2015,Nedic2016}, an algorithm named the distributed dual averaging push method is presented in \cite{Liang2020}, and it is built with the push-sum technique and the dual averaging method.

Many interesting works have been done for distributed optimization problems, and the delayed information needs to be considered in some situations. We find some literatures for the convex optimization problems with delays \cite{Agarwal2012,Wang2015,KI2011,KI2012,Yang2017,Hatanaka2018,DongWang2018,JLi2016,JZhang2016,Doan2017,Assran2020,Zhang2020}. In \cite{Agarwal2012}, authors propose the delayed dual averaging method and the delayed mirror descent method. Then, distributed dual averaging algorithm with subgradient delays is developed in \cite{Wang2015}.
\cite{Doan2017} considers a distributed gradient-based consensus algorithm and analyzes the convergence rate of this algorithm for solving distributed optimization problems over the network with communication delays that are inevitable in distributed systems.
\cite{Assran2020} and \cite{Zhang2020} investigate convergence results of asynchronous distributed convex optimization algorithms which are caused by communication delays in multi-agent systems.

This article investigates the distributed convex optimization problem for the multi-agent system with communication delays over a sequence of time-varying directed graphs. The main contribution can be summarized as the following points.
\begin{enumerate}
\item Firstly,  it is assumed that communication delays exist in investigated multi-agent systems and the fixed delay model introduced in \cite{KI2011} is applied in our work. We find that the maximum value of the communication delay on the edge affects the convergence result of the proposed algorithm.
\item Secondly, our research is studied over multi-agent systems whose components are multiple nodes, and all the node are available to their own objective functions, and can only share information with neighbours over strongly connected and time-varying unbalanced directed networks. For the fact that all the nodes in the system are only able to communicate over the unbalanced network, we apply the push-sum protocol for this situation which can compute weights with merely out degrees of nodes and eliminate the request of doubly stochastic weight matrices.
\item Last but not least, the PS-DDA algorithm which combines the push-sum method and the dual averaging method is employed in this paper. It is proved that PS-DDA algorithm is able to converge at a rate of $\calO\left(T^{-0.5}\right)$ with proper step size, and $T$ denotes iteration span.
\end{enumerate}


The remaining paper is organized as the following: The introduction to PS-DDA algorithm and fixed delay model is illustrated in Section II. The main results under the delay model are in Section III, providing a detailed convergence proof of PS-DDA algorithm. Some comparisons between simulation results of PS-DDA algorithm with and without communication delays are presented in Section IV, discussion of the convergence results of PS-DDA algorithm and other algorithms are proposed as well. In the last but not the least part, Section V gives a conclusion of our research, and proposes some future directions.

{\it Notations:} $\mathbb{R}$ represents the real number set, and $d$ is dimension number. $Q$ is used to represent a matrix, the transpose of this matrix can be written as $Q^{T}$. and the $i$th column of this matrix is $\left[Q\right]_i$. We let $\mathbf{1}$ denote a vector of all ones. For two positive integers $a$ and $b$, $\mathbf{I}_{a}$ denotes the $a\times a$ identity matrix, and $\mathbf{0}_{a\times b}$ is the $a\times b$ zero matrix. $x^T$ is the transpose of vector $x$, whose $i$th element is $\left[x\right]_{i}$. $\left\Vert x\right\Vert_{1}$ is defined as the 1-norm on $\mathbb{R}^d$. $\left\langle z,x\right\rangle=z^{T}y$ is defined as the inner product of two vectors. $\left\Vert x\right\Vert=\left\langle x,x\right\rangle^{1/2} $ is the standard Euclidean norm and $\left\Vert x\right\Vert_{*}=\sup_{\left\Vert y=1\right\Vert }\left\langle x,y\right\rangle $ denotes the dual norm to the standard Euclidean norm. $f_1\left(T\right)=\mathcal{O}\left(f_2\left(T\right)\right)$ is used to express the relationship between two functions $f_1$ and $f_2$ when there is $t<\infty$ and positive constant $C<\infty$ for $T\geq t$ to meet the requirement that $f_1\left(T\right)\leq C f_2\left(T\right)$. For the function $f$, $\nabla f$ denotes its subgradient.

\section{Preliminaries and Formulation}

\subsection{Problem Formulation}
The distributed convex optimization problem over a multi-agent system is considered in our paper. $\mathcal{G}\left(t\right)=(\mathcal{V},\mathcal{E}\left(t\right))$ denotes a time-varying directed graph over the set of vertex $\mathcal{V}=\left\{1,\ldots,m\right\}$ and the set of edge $\mathcal{E}\left(t\right)$. We assume that there are $N$ directed edges in the union network and the set $\bigcup^{(t+1)B-1}_{i=tB} \mathcal{E}\left(i\right)=\left\{\mathcal{E}_1,\ldots,\mathcal{E}_N\right\}$. $\mathcal{N}_{i}^{\mathrm{in}}\left(t\right)=\left\{ j|\left(j,i\right)\in \mathcal{E}\left(t\right)\right\} \cup\left\{ i\right\} $ and $\mathcal{N}_{i}^{\mathrm{out}}\left(t\right)=\left\{ j|\left(i,j\right)\in \mathcal{E}\left(t\right)\right\} \cup\left\{ i\right\} $ denote node $i$'s in-neighbours and out-neighbours at time $t$, respectively. Every agent $i\in\left\{1,\ldots,m\right\}$ in multi-agent system possesses its own objective function: $f_{i}:\mathbb{R}^{d}\rightarrow\mathbb{R}$, and the main purpose is solving:
\BEAS
\min f\left(x\right)=\frac{1}{m}\sum_{i=1}^{m}f_{i}\left(x\right),\;x\in\mathcal{X},
\label{obj-fun}
\EEAS
where the set $\mathcal{X}$ is closed and convex.

The following assumption is proposed for the network graphs of the multi-agent system.
\begin{assumption}\label{assump-G}
With the positive integer $B$, for each $t>0$, the graph is strongly connected with the following edge set
\BEASN
\bigcup^{(t+1)B-1}_{i=tB} \mathcal{E}\left(i\right).
\EEASN
\end{assumption}
Considering the objective function $f_i$, $i=1,\dots,m$ and the optimization problem (\ref{obj-fun}), it is assumed that the following assumptions are satisfied.
\begin{assumption}\label{assump:opt}
Suppose that there invariably exist some solutions $x^{*}$ to problem (\ref{obj-fun}), and $X^{*}$ is the set of optimal solutions $x^{*}$.
\end{assumption}
\begin{assumption}\label{assump:convex}
For node $i\in\left\{1,\ldots,m\right\}$, every function $f_i$ is convex, that is,
\[
f_i\left(y\right)-f_i\left(x\right)\geq\nabla f_{i}\left(x\right)^T\left(y-x\right),\qquad \forall x,y\in\mathbb{R}^{d}.
\]
\end{assumption}

\begin{assumption}\label{assump:Lipschitz}
For node $i\in\left\{1,\ldots,m\right\}$, every objective function $f_i$ is $L$-Lipschitz continuous, that means,
\BEASN
\left|f_{i}\left(x\right)-f_{i}\left(y\right)\right| \leq L\left\| x-y\right\|,\qquad \forall x,y\in\mathbb{R}^{d}.
\EEASN
\end{assumption}

Then for any $x\in\mathcal{X}$ and any subgradient $g_{i}\in\nabla f_{i}\left(x\right)$, we have $\left\Vert g_{i}\right\Vert_* \leq L$. Each agent $i$ in a multi-agent system has its own parameter vector $x_i\in\mathbb{R}^{d}$, can only be available to its own loss function $f_i$, and receive messages from its in-neighbor $j\in\mathcal{N}_{i}^{\mathrm{in}}\left(t\right)$. We define $f\left(x\right)=\frac{1}{m}\sum_{i=1}^{m}f_i\left(x\right)$, and the function $f$ is certainly a convex function.
\subsection{PS-DDA Algorithm}

We will talk about the PS-DDA algorithm in this subsection. Each node in a multi-agent system has a local estimate $x_i\left(t\right)$, a weight $w_i\left(t\right)$ and a dual variable $z_i\left(t\right)$ that maintains the accumulated subgradient at $t$. Then, for any arbitrary $t\geq0$ and $i=1,\ldots,m$, node $i$ updates as:
\BEAS
&&w_{i}\left(t+1\right)=\sum_{j=1}^{m}p_{ij}\left(t\right)w_{j}\left(t\right),\nn\\
&&z_{i}\left(t+1\right)=\sum_{j=1}^{m}p_{ij}\left(t\right)z_{j}\left(t\right)+g_{i}\left(t\right),\nn\\
&&x_{i}\left(t+1\right)=\Pi_{\mathcal{X}}^{\psi}\left(\frac{z_{i}\left(t+1\right)}{w_{i}\left(t+1\right)},\alpha\left(t+1\right)\right),
\label{algorithm-0}
\EEAS

where $g_{i}\left(t\right)$ denotes the subgradient of $f_{i}\left(x\right)$ at $x_{i}\left(t\right)$, $\left\{ \alpha\left(t\right)\right\} _{t=1}^{\infty}$ is a sequence of positive step sizes which is non-increasing, and $P\left(t\right)$ is a weight matrix that consists of the weights $p_{ij}\left(t\right)$ used to compute $w_{i}\left(t+1\right)$ and $z_{i}\left(t+1\right)$. It has that $P\left(t\right)=\left[p_{ij}\left(t\right)\right]_{m\times m}$. For $p_{ij}\left(t\right)$, it is supposed to be expressed as:
\BEAS
p_{ij}\left(t\right)=\begin{cases}
\frac{1}{d_{j}\left(t\right)}, &\qquad j\in\mathcal{N}_{i}^{\mathrm{in}}\left(t\right),\\
0, &\qquad \mathrm{otherwise}.
\end{cases}
\label{algorithm-1}
\EEAS
where $d_{j}\left(t\right)$ represents the amount of edges from node $j$ at time $t$ including itself (\ie{, $d_{j}\left(t\right)=\left|\mathcal{N}_{j}^{\mathrm{out}}\left(t\right)\right|$}). $P\left(t\right)$ is a column stochastic matrix satisfying $\mathbf{1}^{T}P\left(t\right)=\mathbf{1}^{T}$ and represents the structure of the graph $\mathcal{G}\left(t\right)$. For all agents $i=1,2,\ldots,m$, we have $x_{i}\left(0\right)\in \mathbb{R}^d$, $w_i\left(0\right)=1$ and $z_i\left(0\right)=0$. The projection operator $\Pi_{\mathcal{X}}^{\psi}\left(\cdot,\cdot\right)$ is defined in the following,
\[
\Pi_{\mathcal{X}}^{\psi}\left(z,\alpha\right)=\underset{x\in\mathcal{X}}{\arg\min}\left\{ \left\langle z,x\right\rangle +\frac{1}{\alpha}\psi\left(x\right)\right\}
\]
where $\psi\left(x\right):\mathbb{R}^{d}\rightarrow\mathbb{R}$ is a proximal function, and it satisfies the following.
\begin{assumption}\label{assump:psi}
 Suppose $\psi\left(x\right):\mathbb{R}^{d}\rightarrow\mathbb{R}$ is $1$-strongly convex, that is $\psi\left(y\right)\geq\psi\left(x\right)+\nabla\psi\left(x\right)^{T}\left(y-x\right)+\frac{1}{2}\left\Vert y-x\right\Vert ^{2}$, $\forall x,y\in\mathbb{R}^{d}$. It is also supposed $\psi\left(x\right)\geq0$, $\psi\left(0\right)=0$, and for $x^{\star}\in X^{\star}$, the bound $\psi\left(x^{*}\right)\leq R^{2}$ is valid.
\end{assumption}

\begin{remark}
Many distributed convex optimization algorithms are assumed to compute over balanced networks \cite{Doan2019,Yuan2018,Xiong2021,KI2012}, and are possible to utilize doubly stochastic weight matrices. However, algorithms relay on doubly stochastic weight matrices might be dissatisfactory in some real-world applications. In order to work over unbalanced communication graphs and remove the request for doubly stochastic weight matrices, the push-sum method which has been employed in many previous works like \cite{nedic2015,Nedic2016} is applied in algorithm (\ref{algorithm-0}). Compare with the distributed dual averaging algorithm with doubly stochastic matrices in \cite{KI2012}, which needs information of graph sequence and the amount of agents to compute weights, the PS-DDA algorithm in this article can calculate weights with merely the information of the out-degree of each node.
\end{remark}

\begin{remark}
In some previous gradient descent method like the projected consensus algorithm in \cite{Nedic2010,Yuan2021,Liu2021}, the projection step used in the update may bring some technical difficulties. However, the distributed dual averaging algorithm utilizes the proximal function allowing us to deal with the problem with non-Euclidean geometry that can reduce computing complexity in some cases such as high-dimensional problems.
\end{remark}


\begin{remark}
We make a slight improvement to the PS-DDA algorithm introduced in \cite{KI2012}. Our paper discusses the distributed convex optimization problem over a sequence of time-varying and uniformly strongly connected directed graphs instead of a connected fixed directed graph supposed in \cite{KI2012}, which will reduce restriction to application to some extent.
\end{remark}

\subsection{Fixed delay Model}

Referring to the fixed delay model in \cite{KI2011}, for a multi-agent system with $m$ nodes, it is supposed that there are $N$ edges in the edge set $\bigcup^{(t+1)B-1}_{i=tB}\mathcal{E}\left(i\right)$ and when the edge $\mathcal{E}_r=\left(i,j\right), r=1,\dots,N$ in the network is connected at time $t$, each agent $j$ receives a message from the agent $i$ with a delay of $\tau_{\mathcal{E}_r}$ time units. This delay is modeled by adding $\tau_{\mathcal{E}_r}$ delay nodes between agents $i$ and $j$ as the delay, and the maximum value of the delay on one edge is supposed to be $\tau_{max}$. With the positive integer $B$, for a network and each $t>0$, we define $\tau=\sum_{\mathcal{E}_r\in\bigcup^{(t+1)B-1}_{i=tB} \mathcal{E}\left(i\right)}\tau_{\mathcal{E}_r}$ as the amount of delay nodes added in edges of the network. At time $t$, if the edge $\left(i,j\right)$ is connected, the message can be sent from node $i$ to its closest delay node between compute nodes $i$ and $j$, if else, the message can not go through to its closest delay node between compute nodes $i$ and $j$. At all the time, all the edges that start from delay nodes in the augmented network are connected and all the delay nodes can send information to its neighbours through the directed strongly connected graph of the network which consists of $m$ compute nodes and $\tau$ delay nodes (See Figure 1).

We describe the time-varying communication topology graph $\mathcal{G}\left(t\right)$ which has no communication delays by a column stochastic matrix $P\left(t\right)$. We augment the network $\mathcal{G}\left(t\right)$ by adding delay nodes so that it can be presented by another column stochastic matrix $Q\left(t\right)$. The operation to develop the matrix $Q\left(t\right)$ from the matrix $P\left(t\right)$ is achieved by inserting delays on edges. Considering that each edge $\mathcal{E}_r=\left(i,j\right),r=1,\dots,N$ has a delay of $\tau_{\mathcal{E}_r}$, the nodes $d_{\left(i,j\right),1},d_{\left(i,j\right),2},\ldots,d_{\left(i,j\right),\tau_{\mathcal{E}_r}}$ are inserted between nodes $i$ and $j$. A message should go through all the delay nodes while transmitting from node $i$ to $j$.

\begin{figure}[H]
\begin{minipage}[t]{0.3\linewidth}
\centering
\rotatebox{360}{\scalebox{0.16}[0.16]{\includegraphics{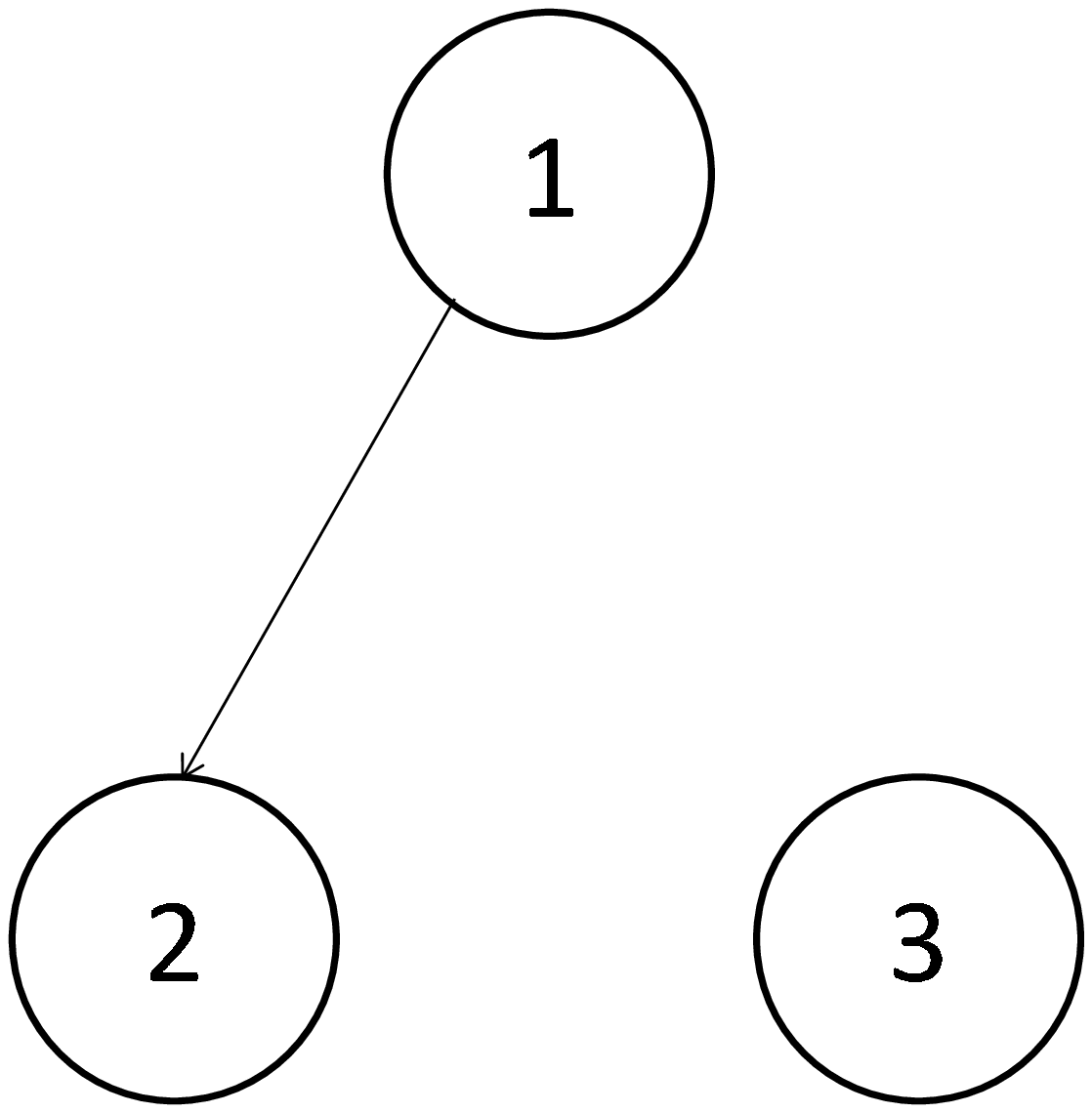}}}
\\[0pt]
\end{minipage}%
\begin{minipage}[t]{0.3\linewidth}
\centering
\rotatebox{360}{\scalebox{0.16}[0.16]{\includegraphics{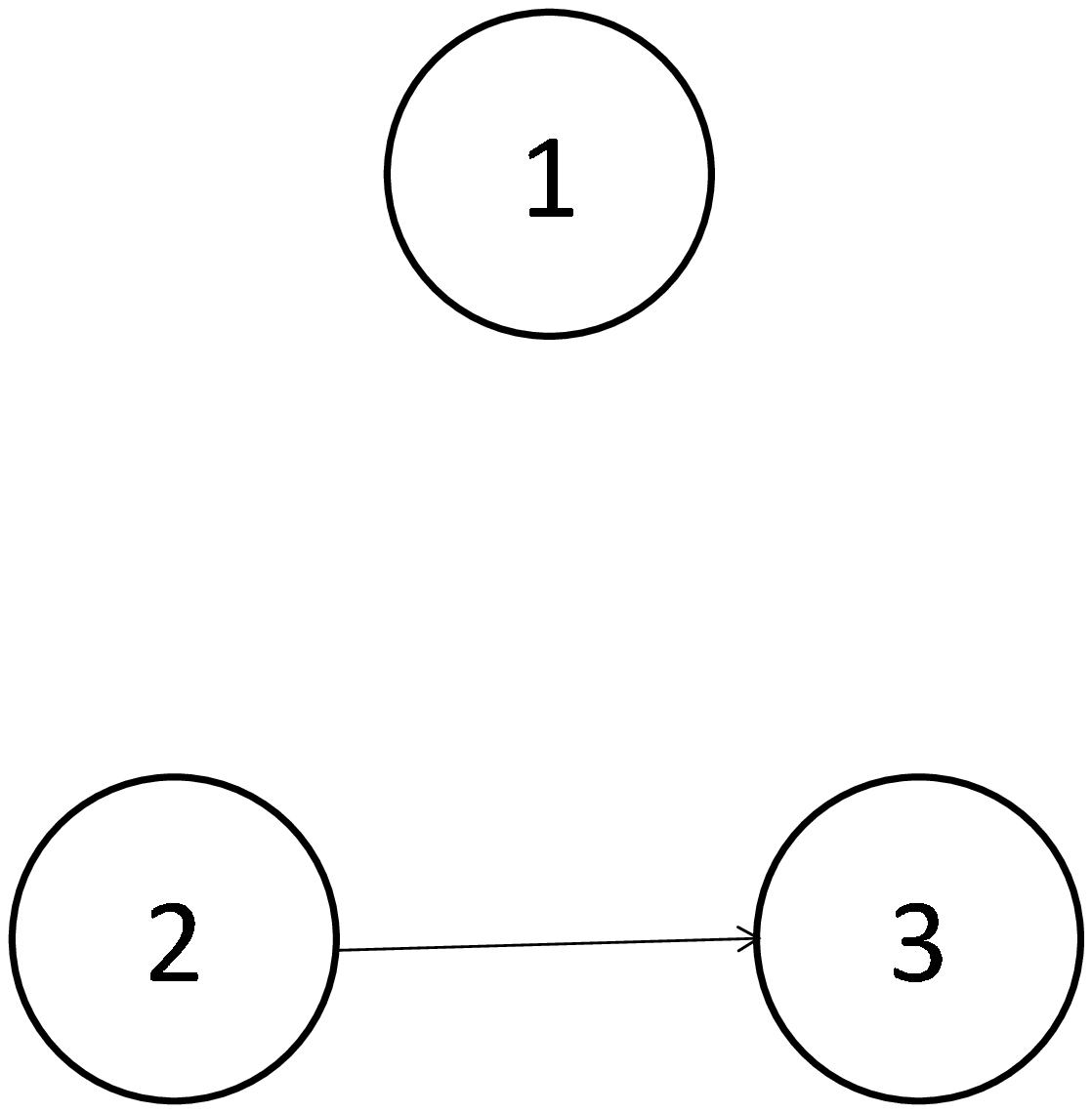}}}
\\[0pt]
\end{minipage}
\begin{minipage}[t]{0.3\linewidth}
\centering
\rotatebox{360}{\scalebox{0.16}[0.16]{\includegraphics{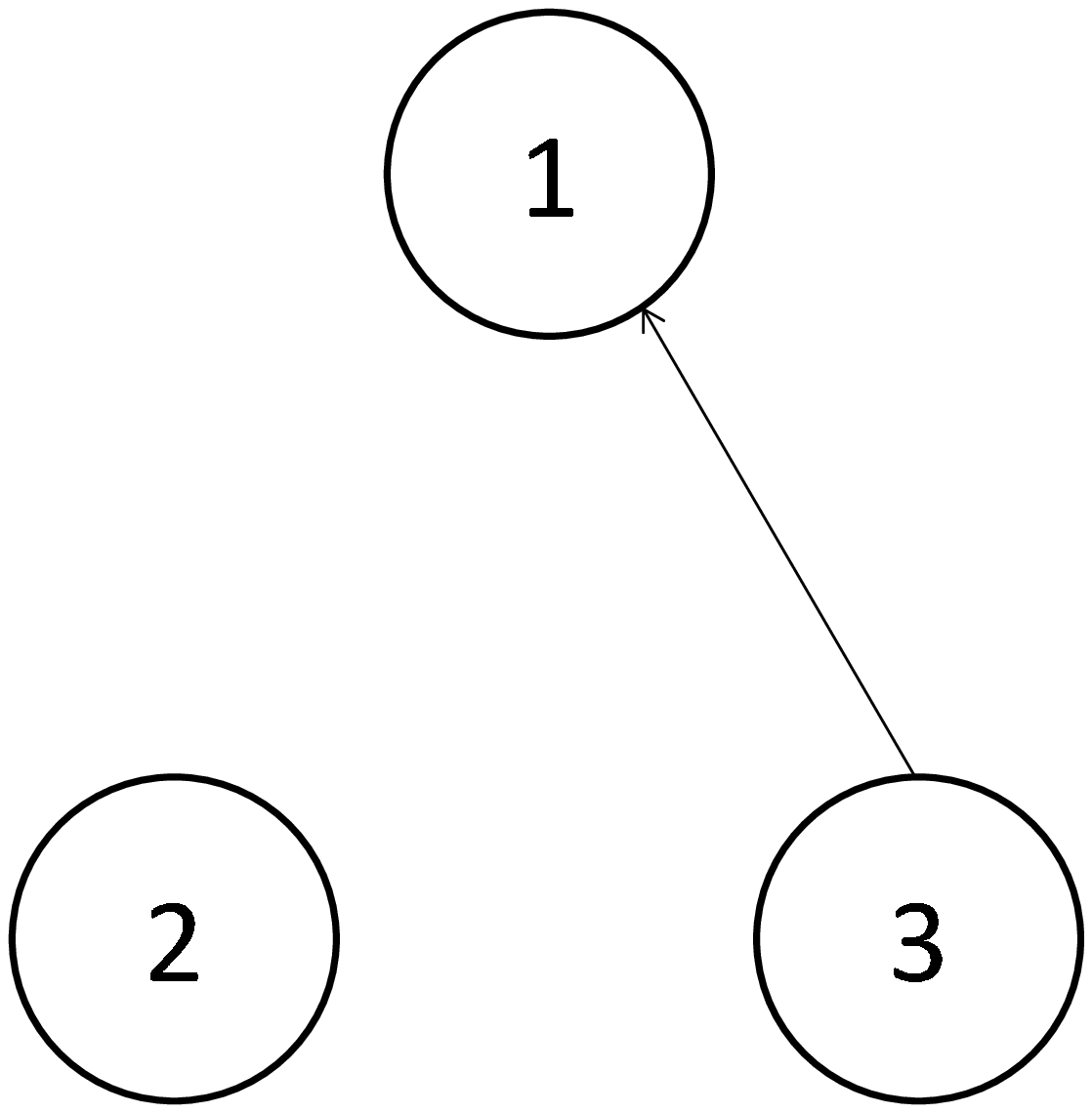}}}
\\[0pt]
\end{minipage}
\centering
\\[3pt]
{\scriptsize{\vc{Figure 1(a). A simple time-varying network with 3 nodes.} }}
\end{figure}

\begin{figure}[H]

\begin{minipage}[t]{0.3\linewidth}
\centering
\rotatebox{360}{\scalebox{0.16}[0.16]{\includegraphics{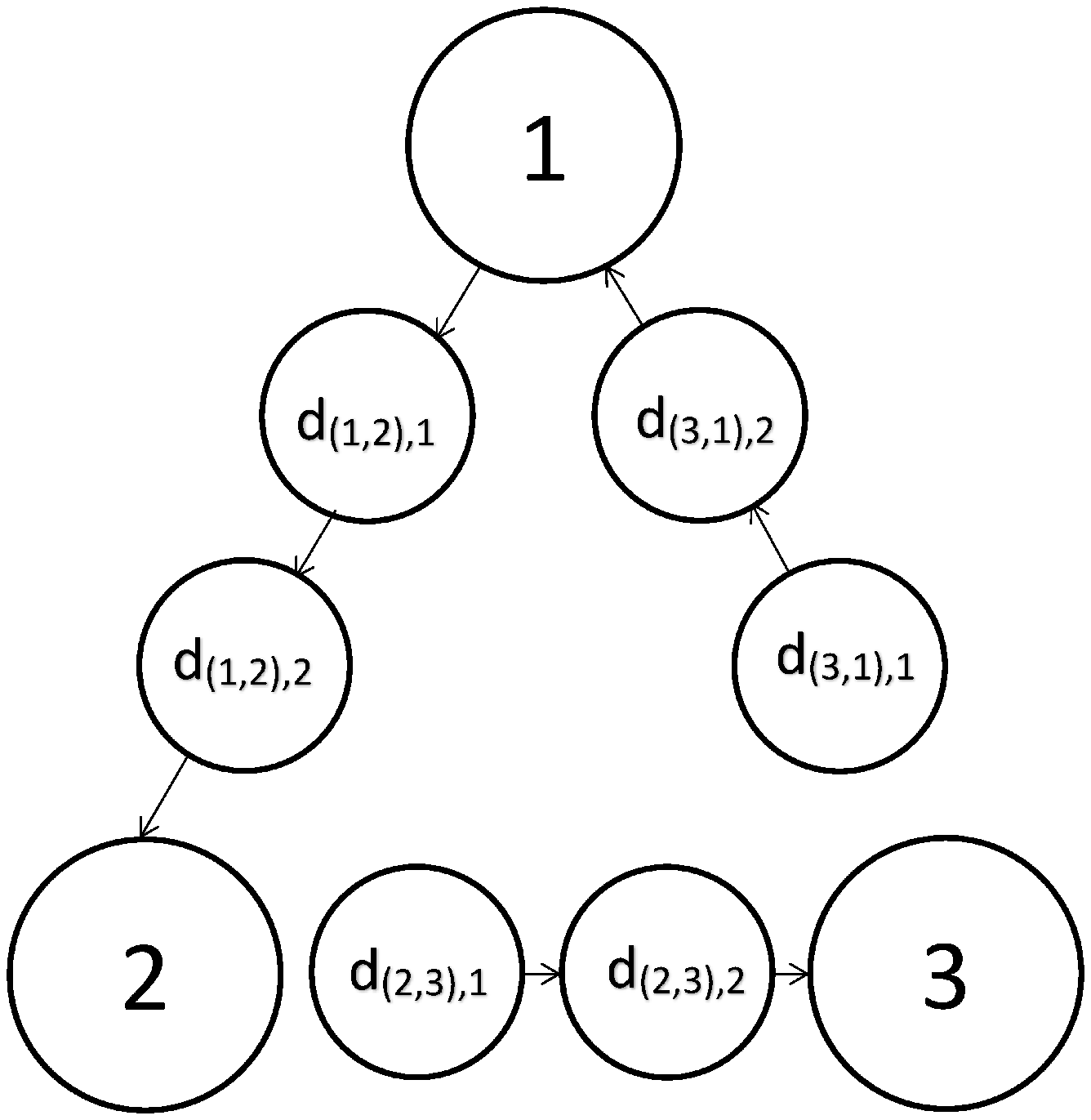}}}
\\[0pt]
\end{minipage}%
\begin{minipage}[t]{0.3\linewidth}
\centering
\rotatebox{360}{\scalebox{0.16}[0.16]{\includegraphics{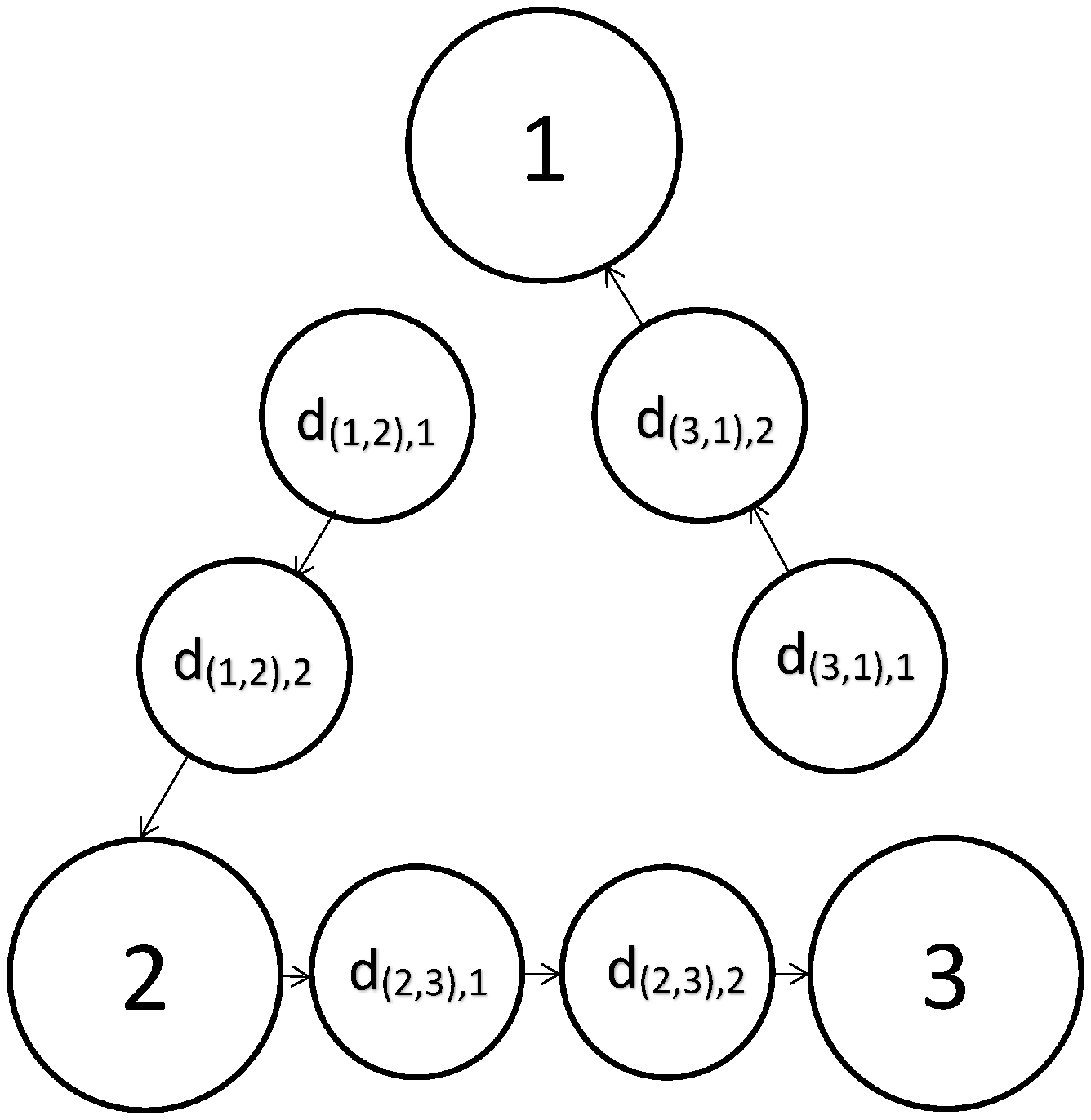}}}
\\[0pt]
\end{minipage}
\begin{minipage}[t]{0.3\linewidth}
\centering
\rotatebox{360}{\scalebox{0.16}[0.16]{\includegraphics{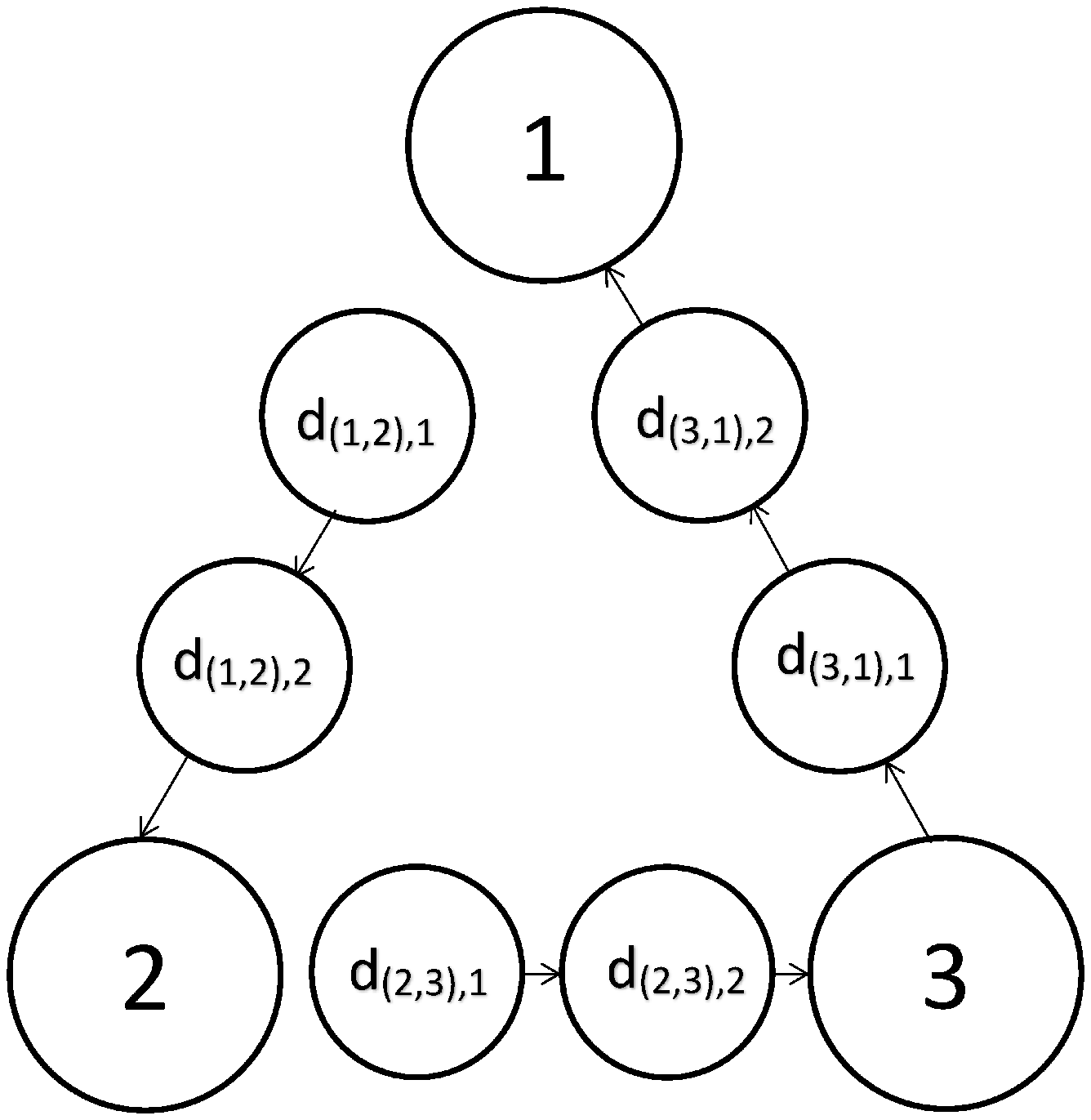}}}
\\[0pt]
\end{minipage}
\centering
\\[3pt]
{\scriptsize{\vc{Figure 1(b). A simple time-varying network with delays.} }}
\end{figure}

The matrix $A\left(t\right)$ is an $m\times m$ square matrix which is used to set the entry of $Q\left(t\right)$ that represents all the directed connections in the network without delay to 0. For $r=1,\ldots,N$, when the edge $\mathcal{E}_r=\left(i,j\right)$ is connected with $\tau_{\mathcal{E}_r}$ delay nodes, then node $i$ passes its information to the closest delay node $d_{\left(i,j\right),1}$. Accordingly, the weight of message sending from node $i$ to delay node $d_{\left(i,j\right),1}$ equals to the one which is applied to deliver message from node $i$ to $j$ directly without any delay, and this procedure is available by a $\tau_{\mathcal{E}_r} \times m$ matrix $C_r\left(t\right)$. When the edge $\mathcal{E}_r=\left(i,j\right)$ is not connected, the $\tau_{\mathcal{E}_r} \times m$ matrix is set to be $\mathbf{0_{\ensuremath{\tau_{\mathcal{E}_{r}}\times m}}}$. The matrix $B_r\left(t\right)$ is an $m \times \tau_{\mathcal{E}_r}$ matrix defined to describe the transmission of information from the last delay node $d_{\tau_{\mathcal{E}_r}}$ to a node $j$. A $\tau_{\mathcal{E}_r}\times\tau_{\mathcal{E}_r}$ square matrix $D_r\left(t\right)$ is utilized to deliver information along the delay chain. $\epsilon_{i}$ represents the $i$-th column of the $m\times m$ identity matrix. $\omega_{\mathcal{E}_r,i}$ denotes the $i$-th column of the $\tau_{\mathcal{E}_r}\times\tau_{\mathcal{E}_r}$ identity matrix. We present the following transformation to achieve the weight matrix $Q\left(t\right)$ which illustrates the graph with $m$ network nodes and $\tau$ delay nodes from the matrix without delay $P\left(t\right)$ at the time $t$.
\BEAS
Q\left(t\right)&=&\left[\begin{array}{c}
\mathbf{I}_{m}\\
\mathbf{0}_{\tau\times m}
\end{array}\right]P\left(t\right)\left[\begin{array}{cc}
\mathbf{I}_{m} & \mathbf{0}_{m\times\tau}\end{array}\right]\nn\\&&+\left[\begin{array}{ccccc}
A\left(t\right) & B_{1}\left(t\right) & B_{2}\left(t\right) & \cdots & B_{N}\left(t\right)\\
C_{1}\left(t\right) & D_{1}\left(t\right) & \mathbf{0} & \cdots & \mathbf{0}\\
C_{2}\left(t\right) & \mathbf{0} & D_{2}\left(t\right) & \vdots & \vdots\\
\vdots & \vdots & \vdots & \ddots & \mathbf{0}\\
C_{N}\left(t\right) & \mathbf{0} & \mathbf{\cdots} & \mathbf{0} & D_{N}\left(t\right)
\end{array}\right],\nn\\
A\left(t\right)&=&-\sum_{\left(i,j\right)\in\mathcal{E}\left(t\right)}\epsilon_{j}\epsilon_{j}^{T}P\left(t\right)\epsilon_{i}\epsilon_{i}^{T},\nn\\
B_{r}\left(t\right)&=&\epsilon_{j}\omega_{\mathcal{E}_{r},\tau_{\mathcal{E}_{r}}}^{T},\quad\forall\mathcal{E}_{r}=\left(i,j\right)\in\bigcup_{i=tB}^{(t+1)B-1}\mathcal{E}\left(i\right),\nn\\
C_{r}\left(t\right)&=&\begin{cases}
\omega_{\mathcal{E}_{r},1}\epsilon_{j}^{T}P\left(t\right)\epsilon_{i}\epsilon_{i}^{T} & \mathcal{E}_{r}=\left(i,j\right)\in\mathcal{E}\left(t\right)\\
\mathbf{0_{\ensuremath{\tau_{\mathcal{E}_{r}}\times m}}} & \mathcal{E}_{r}=\left(i,j\right)\notin\mathcal{E}\left(t\right)
\end{cases},\nn\\
D_{r}\left(t\right)&=&\sum_{k=1}^{{\tau}_{\mathcal{E}_{r}}-1}\omega_{\mathcal{E}_{r},k+1}\omega_{\mathcal{E}_{r},k}^{T}.
\label{transformation-Q}
\EEAS

We then propose a simple example for more intuitive understanding for the transformation (\ref{transformation-Q}).

\begin{example}
We take the simple network with 3 nodes in Figure 1 as an example. There are three directed graphs, denoted $\mathcal{G}_{1}$, $\mathcal{G}_{2}$ and $\mathcal{G}_{3}$ in Figure 1(a). $P_{1}$, $P_{2}$, and $P_{3}$ are used to define weight matrices related to graphs $\mathcal{G}_{1}$, $\mathcal{G}_{2}$ and $\mathcal{G}_{3}$, respectively. It can be seen that the weight matrices $P_1$, $P_2$ and $P_3$ with no delay are written as
\BEASN
P_1=\begin{bmatrix}
\begin{smallmatrix}
0.5 & 0  & 0  \\
0.5 & 1  & 0  \\
0   & 0  & 1
\end{smallmatrix}
\end{bmatrix},
P_2=\begin{bmatrix}
\begin{smallmatrix}
 1   & 0    & 0  \\
 0   & 0.5  & 0  \\
 0   & 0.5  & 1  \\
\end{smallmatrix}
\end{bmatrix},
P_3=\begin{bmatrix}
\begin{smallmatrix}
 1   & 0  & 0.5   \\
 0   & 1  & 0     \\
 0   & 0  & 0.5   \\
\end{smallmatrix}
\end{bmatrix}.
\EEASN

From Figure 1, we notice that in the Figure 1(b), there exist two delay nodes between every edge. When node $i$ gives out a message to node $j$, this message is passed through delay nodes $d_{(i,j),1}$ and $d_{(i,j),2}$ and then received by node $j$. The augmented graph can be represented by column stochastic weight matrices $Q_1$, $Q_2$ and $Q_3$.
\BEASN
&&Q_1=
\begin{bmatrix}
\begin{smallmatrix}
0.5 & 0  & 0   &0           &0           &0           &0           &0           &1\\
0   & 1  & 0   &0           &1           &0           &0           &0           &0\\
0   & 0  & 1   &0           &0           &0           &1           &0           &0\\
0.5 & 0  & 0   &0           &0           &0           &0           &0           &0\\
0   & 0  & 0   &1           &0           &0           &0           &0           &0\\
0   & 0  & 0   &0           &0           &0           &0           &0           &0\\
0   & 0  & 0   &0           &0           &1           &0           &0           &0\\
0   & 0  & 0   &0           &0           &0           &0           &0           &0\\
0   & 0  & 0   &0           &0           &0           &0           &1           &0
\end{smallmatrix}
\end{bmatrix},
Q_2=\begin{bmatrix}
\begin{smallmatrix}
1   & 0    & 0   &0           &0           &0           &0           &0           &1\\
0   & 0.5  & 0   &0           &1           &0           &0           &0           &0\\
0   & 0    & 1   &0           &0           &0           &1           &0           &0\\
0   & 0    & 0   &0           &0           &0           &0           &0           &0\\
0   & 0    & 0   &1           &0           &0           &0           &0           &0\\
0   & 0.5  & 0   &0           &0           &0           &0           &0           &0\\
0   & 0    & 0   &0           &0           &1           &0           &0           &0\\
0   & 0    & 0   &0           &0           &0           &0           &0           &0\\
0   & 0    & 0   &0           &0           &0           &0           &1           &0
\end{smallmatrix}
\end{bmatrix},\nn\\
&&Q_3=\begin{bmatrix}
\begin{smallmatrix}
1   & 0  & 0   &0           &0           &0           &0           &0           &1\\
0   & 1  & 0   &0           &1           &0           &0           &0           &0\\
0   & 0  & 0.5 &0           &0           &0           &1           &0           &0\\
0   & 0  & 0   &0           &0           &0           &0           &0           &0\\
0   & 0  & 0   &1           &0           &0           &0           &0           &0\\
0   & 0  & 0   &0           &0           &0           &0           &0           &0\\
0   & 0  & 0   &0           &0           &1           &0           &0           &0\\
0   & 0  & 0.5 &0           &0           &0           &0           &0           &0\\
0   & 0  & 0   &0           &0           &0           &0           &1           &0
\end{smallmatrix}
\end{bmatrix}.
\EEASN
\end{example}

\begin{remark}
With the fixed delay model, it provides an idea to tackle the communication delays in the convex optimization problems by adding delay nodes between compute nodes. It is notable that we suggest that the self loop message does not have delay, which means each node is invariably available to its latest local estimation. Using this way to deal with communication delays has a shortcoming that if the network has very long communication delays, the number of calculations will be quite large for the reason that all the compute nodes and delay nodes in the network have to compute for every time of iteration.
\end{remark}


We use $\tau$ delay nodes to model communication delays in the network. It is supposed that every delay node has a cost function $f_{i}\left(x\right)=0$, $i=m+1,\ldots,m+\tau$. For $i\in\left\{1,\ldots,m\right\}$, $g_{i}\left(t\right)$ is defined as the subgradient of $f_{i}\left(x\right)$ at $x_{i}\left(t\right)$; For $i\in\left\{m+1,\ldots,m+\tau\right\}$, $g_{i}\left(t\right)=0$ is supposed to be valid. For $i\in\left\{1,\ldots,m\right\}$, the initial values $w_i\left(0\right)=1$ and $z_i\left(0\right)=0$; for $i\in\left\{m+1,\ldots,m+\tau\right\}$, the initial values $w_i\left(0\right)=0$ and $z_i\left(0\right)=0$. Matrix $Q\left(t\right)=\left[q_{ij}\left(t\right)\right]_{(m+\tau)\times (m+\tau)}$ is used as a transition matrix in PS-DDA algorithm with communication delays, replacing matrix $P\left(t\right)$ in algorithm (\ref{algorithm-0}). Then, for any arbitrary $t\geq0$ and $i\in\left\{1,\dots,m,m+1,\ldots,m+\tau\right\}$, the weight $w_i\left(t\right)$ and the dual variable $z_i\left(t\right)$ are updated as,
\BEAS
&&w_{i}\left(t+1\right)=\sum_{j=1}^{m+\tau}q_{ij}\left(t\right)w_{j}\left(t\right),\nn\\
&&z_{i}\left(t+1\right)=\sum_{j=1}^{m+\tau}q_{ij}\left(t\right)z_{j}\left(t\right)+g_{i}(t),
\label{algorithm-2}
\EEAS
For any arbitrary $t\geq0$ and $i\in\left\{1,\ldots,m\right\}$, the local estimate $x_{i}\left(t+1\right)$ is generated by
\BEAS
x_{i}\left(t+1\right)=\Pi_{\mathcal{X}}^{\psi}\left(\frac{z_{i}\left(t+1\right)}{w_{i}\left(t+1\right)},\alpha\left(t+1\right)\right).
\label{algorithm-3}
\EEAS

We then define the state transition $Q\left(t:s\right)\in\mathbb{R}^{\left(m+\tau\right)\times\left(m+\tau\right)}$ as
\BEASN
Q\left(t:s\right)=\begin{cases}
Q\left(t\right)Q\left(t-1\right)\cdots Q\left(s\right), & t\geq s\\
\mathbf{I}_{\left(m+\tau\right)}, & t<s
\end{cases}
\EEASN


\begin{remark}
In practical applications, it is usual to have communication delays, which affect computation performances of distributed convex optimization algorithms. The fixed delay model is a necessary method to deal with communication delays between nodes, and have been applied in many works \cite{KI2012,Assran2020,Zhang2020}. The fixed delay model can convert a complicated system with time delays which may receive messages from different time to one without communication delays. Therefore, the fixed delay model reduces the difficulty of computation of the delayed system by increasing the number of nodes in the network.
\end{remark}

\section{The Convergence}
In this section, we are going to introduce our research achievement which is proposed by theorems. These theorems reflect the convergence result of the proposed algorithm while considering communication delays. Our results and processes of their proofs will be presented in detail in the next two subsections.
\subsection{Main Results}

We first define the local running average as $\hat{x}_{i}\left(T\right)=\frac{1}{T}\sum_{t=1}^{T}x_{i}\left(t\right)$ and the average cumulated gradient as $\bar{z}\left(t\right)=\frac{1}{m}\sum_{i=1}^{m+\tau}z_{i}\left(t\right)$, then a theorem possessing that the cost function value of the local running average $\hat{x}_{i}\left(T\right)$ converges to the optimum with the PS-DDA algorithm with communication delays.
\begin{theorem}\label{theorem-basic}
Suppose Assumptions \ref{assump-G}-\ref{assump:Lipschitz} are established, and the sequences $\left\{z_i\left(t\right)\right\}_{t=1}^\infty$ and $\left\{x_i\left(t\right)\right\}_{t=1}^{\infty}$ are generated by the algorithms (\ref{algorithm-2})-(\ref{algorithm-3}), with the non-increasing step size sequence $\left\{\alpha\left(t\right)\right\}_{t=1}^\infty$. For every node $i\in \left\{1,2,\ldots,m\right\}$ and $T\geq1$,
\BEAS
&&f\left(\hat{x}_{i}\left(T\right)\right)-f\left(x^{*}\right)\nn\\
&\leq&\frac{L}{T}\sum_{t=1}^{T}\alpha\left(t\right)\left\Vert \frac{z_{i}\left(t\right)}{w_{i}\left(t\right)}-\bar{z}\left(t\right)\right\Vert _{*}\nn\\
&&+\frac{L^{2}}{2T}\sum_{t=1}^{T}\alpha\left(t\right)+\frac{1}{T\alpha\left(T\right)}\psi\left(x^{*}\right)\nn\\
&&+\frac{2L}{mT}\sum_{t=1}^{T}\sum_{i=1}^{m}\alpha\left(t\right)\left\Vert \frac{z_{j}\left(t\right)}{w_{j}\left(t\right)}-\bar{z}\left(t\right)\right\Vert _{*}.
\label{theorem-basic-0}
\EEAS
\end{theorem}

In the result of above Theorem, the bound of $f\left(\hat{x}_{i}\left(T\right)\right)-f\left(x^{*}\right)$ is related to the network error term $\left\Vert \frac{z_{i}\left(t\right)}{w_{i}\left(t\right)}-\bar{z}\left(t\right)\right\Vert_{*} $. Then we try to bound the network error $\left\Vert \frac{z_{i}\left(t\right)}{w_{i}\left(t\right)}-\bar{z}\left(t\right)\right\Vert_{*} $ and derive convergence rates which have a relationship with the network characteristics in the following theorem.

\begin{theorem}\label{theorem-alpha}
Suppose Assumptions \ref{assump-G}-\ref{assump:psi} are established, and the sequences $\left\{z_i\left(t\right)\right\}_{t=1}^\infty$ and $\left\{x_i\left(t\right)\right\}_{t=1}^\infty$ are generated by algorithms (\ref{algorithm-2}) and (\ref{algorithm-3}) respectively, with step sizes $\left\{\alpha\left(t\right)\right\}_{t=1}^\infty$ which satisfy $\alpha\left(t\right)=\infty$, $\sum_{t=1}^{\infty}\alpha\left(t\right)^{2}<\infty$, and $\alpha\left(t\right)\leq\alpha\left(s\right)$, $\forall t\geq s\geq1$. Then for every node $i\in \left\{1,2,\ldots,m\right\}$ and $T\geq1$,
\BEAS
&&f\left(\hat{x}_{i}\left(T\right)\right)-f\left(x^{*}\right)\nn\\
&\leq&\frac{L^{2}}{2T}\sum_{t=1}^{T}\alpha\left(t\right)+\frac{R^{2}}{T\alpha\left(T\right)}+\frac{3\varGamma L^{2}}{T}\sum_{t=1}^{T}\alpha\left(t\right),
\label{theorem-alpha-0}
\EEAS
where $\Gamma=\frac{mC}{\delta}\left(\frac{1}{\left(1-\lambda\right)\lambda}+t^{*}\lambda^{t^{*}-1}\right)$, here $C=4\frac{1+\left(\frac{1}{m}\right)^{-\Omega}}{1-\left(\frac{1}{m}\right)^{\Omega}}$, $\lambda=\left(1-\frac{1}{m^{\Omega}}\right)^{\frac{1}{\Omega}}$, $\delta\geq1/m^{\Omega+1}$, $\Omega=\left(m-1\right)B+m\left(\tau_{max}+1\right) \geq 3$ and $t^{*}=\arg\max\left(t\lambda^{t-1}\right)$.
\end{theorem}

\begin{remark}
From the result of Theorem \ref{theorem-alpha}, it is obvious that when the maximum value of the delay on one edge $\tau_{max}$ is lager, then $\Omega=\left(m-1\right)B+m\left(\tau_{max}+1\right)$ will be lager. Furthermore, the value of $\Omega$ affects the value of $C$, $\lambda$ and the lower bound of $\delta$. If the value of $\Omega$ increases, the lower bound of $\delta$ which is $1/m^{\Omega+1}$ will decrease. When we consider $\Omega\geq3$ as a variable, $C=4\frac{1+\left(\frac{1}{m}\right)^{-\Omega}}{1-\left(\frac{1}{m}\right)^{\Omega}}$ and $\lambda=\left(1-\frac{1}{m^{\Omega}}\right)^{\frac{1}{\Omega}}$ as two functions related to $\Omega$, we can find that both values of functions $C$ and $\lambda$ grow as variable $\Omega$ increases. With the definition of function $\Gamma=\frac{mC}{\delta}\left(\frac{1}{\left(1-\lambda\right)\lambda}+t^{*}\lambda^{t^{*}-1}\right)$, which can be regarded as a function related to $\lambda$, it can be computed that $\Gamma$ is going to rise as $\lambda$ increases. Therefore, it is concluded that the maximum delay $\tau_{max}$ influences the convergence result of the proposed algorithm. If the maximum delay $\tau_{max}$ becomes larger, the upper bound of $f\left(\hat{x}_{i}\left(T\right)\right)-f\left(x^{*}\right)$ will get larger, which means the convergence result of the algorithm is going to be worse.
\end{remark}

In Theorem \ref{theorem-alpha}, if we choose proper step size $\alpha\left(t\right)$, the bound of $f\left(\hat{x}_{i}\left(T\right)\right)-f\left(x^{*}\right)$ will present the best performance and the algorithm converge with a rate of $\calO\left(T^{-0.5}\right)$. To achieve the best bound, we first set $\alpha\left(t\right)=\frac{\Lambda}{\sqrt{t}}$ for $t\geq1$ where $\Lambda>0$ is a constant, and with the fact that $\sum_{t=1}^T \frac{1}{\sqrt{t}}\leq 2\sqrt{t}$, we compute
\[
\frac{L^{2}}{2T}\sum_{t=1}^{T}\alpha\left(t\right)+\frac{3\varGamma L^{2}}{T}\sum_{t=1}^{T}\alpha\left(t\right)=\frac{R^{2}}{T\alpha\left(T\right)}.
\]
Then the optimal selection of $\alpha\left(t\right)$ can be obtained that,
\[
\alpha\left(t\right)=\frac{R}{L\left(\sqrt{1+6\Gamma}\right)}\frac{1}{\sqrt{t}}.
\]
As a result, the following Corollary which relates to Theorem \ref{theorem-alpha} holds.
\begin{corollary}
We set that the conditions of Theorem \ref{theorem-alpha} hold and $\alpha\left(t\right)=\frac{R}{L\left(\sqrt{1+6\Gamma}\right)}\frac{1}{\sqrt{t}}$, then for every node $i\in \left\{1,2,\ldots,m\right\}$,
\BEASN
f\left(\hat{x}_{i}\left(T\right)\right)-f\left(x^{*}\right)\leq2RL\sqrt{1+6\Gamma}\frac{1}{\sqrt{T}}.
\EEASN
\end{corollary}

\subsection{Proof of Results}

In the begin part of this section, we would like to demonstrate some lemmas which are solid under the conditions of our paper, and they will be applied to prove the main results in our paper, their proof can be seen in \cite{Duchi2012} and \cite{nedic2010}.
\begin{lemma}\label{lemma-gxxstar}(\cite{Duchi2012})
We set $\left\{ g\left(t\right)\right\} _{t=0}^{\infty}\subset\mathbb{R}^{d}$ as a sequence of vectors, the sequence $\left\{ \alpha\left(t\right)\right\} _{t=1}^{\infty}$ is non-increasing, then for $x^*\in \mathcal{X}$ and following sequence,
\BEAS
\tilde{x}\left(t\right)=\Pi_{\mathcal{X}}^{\psi}\left(\sum_{r=0}^{t}g\left(r\right),\alpha\left(t\right)\right),
\label{lemma-gxxstar-0}
\EEAS
it obtains that
\BEAS
&&\sum_{t=1}^{T}\left\langle g\left(t\right),\tilde{x}\left(t\right)-x^{*}\right\rangle \nn\\
&\leq&\frac{1}{2}\sum_{t=1}^{T}\alpha\left(t\right)\left\Vert g\left(t\right)\right\Vert _{*}^{2}+\frac{1}{\alpha\left(T\right)}\psi\left(x^{*}\right)
\label{lemma-gxxstar-0}
\EEAS
\end{lemma}

\begin{lemma}\label{lemma-uv}(\cite{Duchi2012})
For a couple of vectors $u,v\in \mathbb{R}^{d}$, we have
\BEAS
\left\Vert \Pi_{\mathcal{X}}^{\psi}\left(u,\alpha\right)-\Pi_{\mathcal{X}}^{\psi}\left(v,\alpha\right)\right\Vert \leq\alpha\left\Vert u-v\right\Vert _{*}
\label{lemma-uv-0}
\EEAS
\end{lemma}

\begin{lemma}\label{lemma-transition}(\cite{nedic2010})
Suppose Assumption \ref{assump-G} holds, and a sequence of stochastic vectors is $\tilde{\phi}\left(s\right)=\left(\tilde{\phi}_{1}\left(s\right),\ldots,\tilde{\phi}_{m}\left(s\right),\tilde{\phi}_{m+1}\left(s\right),\ldots,\tilde{\phi}_{m+\tau}\left(s\right)\right)^{T}\in\mathbb{R}^{m+\tau}$, which satisfies $\tilde{\phi}_{j}\left(s\right)\geq0$ and $\sum_{j=1}^{m+\tau}\tilde{\phi}_{j}\left(s\right)=1$, then for arbitrary $i,j\in\left\{1,\ldots,m,m+1,\ldots,m+\tau\right\}$ and $t\geq s \geq 0$, we have
\BEASN
&&\left|\left[Q^{T}\left(t\right)\cdots Q^{T}\left(s\right)\right]_{ij}-\tilde{\phi}_{j}\left(s\right)\right|\nn\\
&\leq&2\frac{1+\left(\frac{1}{m}\right)^{-\Omega}}{1-\left(\frac{1}{m}\right)^{\Omega}}\left(1-\left(\frac{1}{m}\right)^{\Omega}\right)^{\frac{t-s}{\Omega}},
\EEASN
where $\Omega=\left(m-1\right)B+m\left(\tau_{max}+1\right) \geq 3$.
\end{lemma}

Then, we are going to apply Lemma \ref{lemma-transition} to a sequence of graphs that merely satisfy uniformly strongly connected after removing the first several graphs. As a result, the following corollary is proposed whose proof is obvious.
\begin{corollary}\label{corollary-transition}
It is assumed that there is a sequence of graphs $\left\{\mathcal{G}\left(t\right)\right\}$, meeting the requirement: there is an integer $M>0$ such that given any integer $s\geq 0$, a time $0\leq t_{s}\leq M$ for which the graph sequence $\mathcal{G}\left(s+t_s\right)$, $\mathcal{G}\left(s+t_s+1\right)$, $\mathcal{G}\left(s+t_s+2\right)$, $\ldots$ is uniformly strongly connected is existed. Then for each integer $s\geq 0$, we have a sequence of stochastic vectors $\tilde{\phi}\left(s\right)=\left(\tilde{\phi}_{1}\left(s\right),\ldots,\tilde{\phi}_{m}\left(s\right),\tilde{\phi}_{m+1}\left(s\right),\ldots,\tilde{\phi}_{m+\tau}\left(s\right)\right)^{T}\in\mathbb{R}^{m+\tau}$, which satisfies $\tilde{\phi}_{j}\left(s\right)\geq0$ and $\sum_{j=1}^{m+\tau}\tilde{\phi}_{j}\left(s\right)=1$, for arbitrary $i,j\in\left\{1,\ldots,m,m+1,\ldots,m+\tau\right\}$ and $t\geq s \geq 0$,
\BEASN
&&\left|\left[Q^{T}\left(t\right)\cdots Q^{T}\left(s\right)\right]_{ij}-\tilde{\phi}_{j}\left(s\right)\right|\nn\\
&\leq&2\frac{1+\left(\frac{1}{m}\right)^{-\Omega}}{1-\left(\frac{1}{m}\right)^{\Omega}}\left(1-\left(\frac{1}{m}\right)^{\Omega}\right)^{\frac{t-s-M}{\Omega}},
\EEASN
where $\Omega=\left(m-1\right)B+m\left(\tau_{max}+1\right) \geq 3$.
\end{corollary}

\begin{lemma}\label{lemma-lowbounddelta}
For a uniformly strongly connected graph sequence $\left\{\mathcal{G}\left(t\right)\right\}$, we denote
\BEASN
\delta'\triangleq\inf_{t=0,1,\cdots}\left(\min_{1\leq j\leq m}\sum_{i=1}^{m}\left[Q^{T}\left(t\right)\cdots Q^{T}\left(0\right)\right]_{ij}\right),
\EEASN
then we have $\delta'\geq\frac{1}{m^{\Omega+1}}$.
\end{lemma}
{\it Proof:}
From the definition of the weight matrix $Q\left(t\right)$, it is obvious that $\left[Q^{T}\left(t\right)\right]_{ii}=1/d_{i}\left(t\right)\geq\frac{1}{m}$, for all $t$ and $i\in\left\{1,\ldots,m\right\}$. Thus, for $i\in\left\{1,\ldots,m\right\}$ and $t\geq0$,
\BEASN
\left[Q^{T}\left(t+1\right)\cdots Q^{T}\left(0\right)\right]_{ii}\geq\frac{1}{m}\left[Q^{T}\left(t\right)\cdots Q^{T}\left(0\right)\right]_{ii}.
\EEASN
Therefore, it follows that $\sum_{i=1}^{m}\left[Q^{T}\left(t\right)\cdots Q^{T}\left(0\right)\right]_{ij}\geq\frac{1}{m^{\Omega+1}}$ for $j\in\left\{1,\ldots,m\right\}$ and $0\leq t\leq\Omega$. It has been presented in \cite{Jadbabaie2003} that for $t>\Omega$ and $i,j\in\left\{1,\ldots,m\right\}$, $\left[Q^{T}\left(t\right)\cdots Q^{T}\left(0\right)\right]_{ij}$ is positive and has value at least $\frac{1}{m^{\Omega+1}}$. Thus, it holds the bound $\delta'\geq\frac{1}{m^{\Omega+1}}$. The proof is complete.
\hfill$\blacksquare$

Then, applying Corollary \ref{corollary-transition} and Lemma \ref{lemma-lowbounddelta}, the following Corollary is obtained.

\begin{corollary}\label{Collary-bound}
\hspace*{\fill}
\begin{enumerate}
  \item It is supposed that Assumption \ref{assump-G} holds, and a sequence of stochastic vector is defined as $\phi\left(t\right)=\left(\phi_{1}\left(t\right),\ldots,\phi_{m}\left(t\right),\phi_{m+1}\left(t\right),\ldots,\phi_{m+\tau}\left(t\right)\right)^{T}\in\mathbb{R}^{m+\tau}$, which satisfies $\phi_{i}\left(t\right)\geq0$ and $\sum_{i=1}^{m+\tau}\phi_{i}\left(t\right)=1$, then for arbitrary $i,j\in\left\{1,\ldots,m,m+1,\ldots,m+\tau\right\}$, it is obtained that,
      \BEAS
      \left|\left[Q\left(t:s\right)\right]_{ij}-\phi_{i}\left(t\right)\right|\leq C\lambda^{t-s},
      \label{Qts-bound}
      \EEAS
      where $C=4\frac{1+\left(\frac{1}{m}\right)^{-\Omega}}{1-\left(\frac{1}{m}\right)^{\Omega}}$ and $\lambda=\left(1-\frac{1}{m^{\Omega}}\right)^{\frac{1}{\Omega}}$.
  \item We define the quantity $\delta\triangleq\inf_{t=0,1,\cdots}\left(\min_{1\leq i\leq m}\sum_{j=1}^{m}\left[Q\left(t\right)\cdots Q\left(0\right)\right]_{ij}\right)$, which satisfies $\delta\geq\frac{1}{m^{\Omega+1}}$, and it implies that $\sum_{j=1}^{m}\left[Q\left(t:0\right)\right]_{ij}\geq\delta$ for all the $t\geq0$ and $i\in\left\{1,\ldots,m\right\}$.
\end{enumerate}
\end{corollary}
{\it Proof:}
When $t<s$, the inequality (\ref{Qts-bound}) is clearly established; and when $t\geq s\geq 0$, the matrix $A\left(t-a\right)$ is defined as
\BEASN
A\left(t-a\right)=Q\left(s+a\right),\;0\leq a\leq t-s,
\EEASN
and $\phi\left(t\right)=\tilde{\phi}\left(s\right)$. Then for $t\geq s\geq 0$,
\BEAS
Q\left(t:s\right)&=&Q\left(t\right)\cdots Q\left(s\right)=\left[Q^{T}\left(s\right)\cdots Q^{T}\left(t\right)\right]^{T}\nn\\
&=&\left[A^{T}\left(t\right)\cdots A^{T}\left(s\right)\right]^{T},
\label{trans-Q-ts}
\EEAS
it can be seen that the order in which the matrices appear has been opposite. The orientation of each edge in every graph $\mathcal{G}\left(t\right)$ has been effectively reversed by taking the transposes of each matrix. In case of $B=1$, it is obvious that the graph is connected at any time, consequently Lemma \ref{lemma-transition} can be applied. Otherwise, the reversed sequence can be strongly connected if at most the first $B$ graphs of this sequence are removed, then applying Corollary \ref{corollary-transition}, it is obvious that
\BEASN
&&\left|\left[Q\left(t:s\right)\right]_{ij}-\phi_{i}\left(t\right)\right|\nn\\
&=&\left|\left[\left[A^{T}\left(t\right)\cdots A^{T}\left(s\right)\right]^{T}\right]_{ij}-\tilde{\phi}_{i}\left(s\right)\right|\nn\\
&=&\left|\left[A^{T}\left(t\right)\cdots A^{T}\left(s\right)\right]_{ji}-\tilde{\phi}_{i}\left(s\right)\right|\nn\\
&\leq&2\frac{1+\left(\frac{1}{m}\right)^{-\Omega}}{1-\left(\frac{1}{m}\right)^{\Omega}}\left(1-\left(\frac{1}{m}\right)^{\Omega}\right)^{\frac{t-s-B}{\Omega}}\nn\\
&\leq&4\frac{1+\left(\frac{1}{m}\right)^{-\Omega}}{1-\left(\frac{1}{m}\right)^{\Omega}}\left(1-\left(\frac{1}{m}\right)^{\Omega}\right)^{\frac{t-s}{\Omega}},
\EEASN
the first item in Corollary \ref{Collary-bound} is completed and we are going to prove the second item. Using the definition of $\delta$ and Lemma \ref{lemma-lowbounddelta},
\BEASN
\delta
=\inf_{t=0,1,\cdots}\left(\min_{1\leq j\leq m}\sum_{i=1}^{m}\left[A^{T}\left(t\right)\cdots A^{T}\left(0\right)\right]_{ij}\right)\geq\frac{1}{m^{\Omega+1}}.
\EEASN
With the obvious inequality $\sum_{j=1}^{m}\left[Q\left(t:0\right)\right]_{ij}\geq \delta$, the proof is complete.
\hfill$\blacksquare$

Before proving the main results in our paper, we propose a definition to an auxiliary sequence $y\left(t\right)$,
\BEAS
y\left(t\right)=\Pi_{\mathcal{X}}^{\psi}\left(\bar{z}\left(t\right),\alpha\left(t\right)\right).
\label{def-y}
\EEAS
$\bar{z}\left(t\right)=\left(1/m\right)\sum_{i=1}^{m+\tau}z_{i}\left(t\right)$ is defined in previous, and it can be represented as,
\BEAS
\bar{z}\left(t+1\right)
&=&\frac{1}{m}\sum_{j=1}^{m+\tau}\left(\sum_{i=1}^{m+\tau}q_{ij}\left(t\right)\right)z_{j}\left(t\right)+\frac{1}{m}\sum_{i=1}^{m+\tau}g_{i}\left(t\right).\nn\\
\
\label{z-average-0}
\EEAS
The matrix $Q\left(t\right)$ is column stochastic, then $\mathbf{1}^{T}Q\left(t\right)=\mathbf{1}^{T}$ holds which depicts that $\sum_{i=1}^{m+\tau}q_{ij}\left(t\right)=1$. With this function, we obtain
\BEAS
\bar{z}\left(t+1\right)=\frac{1}{m}\sum_{j=1}^{m+\tau}z_{j}\left(t\right)+\frac{1}{m}\sum_{i=1}^{m+\tau}g_{i}\left(t\right)
\label{z-average-1}
\EEAS
Obviously, $g_i\left(t\right)=0$ for $i\in\left\{m+1,\ldots,m+\tau\right\}$, and using update of average cumulated gradient $\bar{z}\left(t\right)$ in equation (\ref{z-average-1}) with the fact that $z_i\left(0\right)=0$, thus for $t\geq0$, $\bar{z}\left(t\right)$ and $y\left(t\right)$ could be rewritten as,
\BEAS
&&\bar{z}\left(t+1\right)=\frac{1}{m}\sum_{s=0}^{t}\sum_{i=1}^{m+\tau}g_{i}\left(s\right),\nn\\
&&y\left(t\right)=\Pi_{\mathcal{X}}^{\psi}\left(\frac{1}{m}\sum_{s=0}^{t-1}\sum_{i=1}^{m+\tau}g_{i}\left(s\right),\alpha\left(t\right)\right).
\label{z-average-y}
\EEAS

After the above preparation, we then come to prove Theorem \ref{theorem-basic}.

{\it Proof of Theorem \ref{theorem-basic}:}

As supposed in our paper, the function $f\left(x\right)$ is Lipschitz continuous and convex, for $T\geq1$
\BEAS
f\left(\hat{x}_{i}\left(T\right)\right)-f\left(x^{*}\right)
&\leq&\frac{1}{T}\sum_{t=1}^{T}\left[f\left(y\left(t\right)\right)-f\left(x^{*}\right)\right]\nn\\
&&+\frac{L}{T}\sum_{t=1}^{T}\left\Vert x_{i}\left(t\right)-y\left(t\right)\right\Vert.
\label{theorem-basic-1}
\EEAS
To bound $\frac{1}{T}\sum_{t=1}^{T}\left[f\left(y\left(t\right)\right)-f\left(x^{*}\right)\right]$, the term $\frac{1}{T}\sum_{t=1}^{T}\frac{1}{m}\sum_{i=1}^{m}f_{i}\left(x_{j}\left(t\right)\right)$ is added and substracted at the same time, then
\BEAS
&&\frac{1}{T}\sum_{t=1}^{T}\left[f\left(y\left(t\right)\right)-f\left(x^{*}\right)\right]\nn\\
&\leq&\frac{1}{T}\sum_{t=1}^{T}\frac{1}{m}\sum_{i=1}^{m}\left[f_{i}\left(x_{j}\left(t\right)\right)-f_{i}\left(x^{*}\right)\right]\nn\\
&&+\frac{1}{T}\sum_{t=1}^{T}\frac{1}{m}\sum_{i=1}^{m}\left[f_{i}\left(y\left(t\right)\right)-f_{i}\left(x_{j}\left(t\right)\right)\right].
\label{theorem-basic-2}
\EEAS
With the convex property of each cost function $f_{i}\left(x\right)$ and $g_{i}\in\nabla f_{i}\left(x_{i}\left(t\right)\right)$ for $i\in\left\{1,\ldots,m\right\}$,
\BEAS
&&\frac{1}{T}\sum_{t=1}^{T}\frac{1}{m}\sum_{i=1}^{m}\left[f_{i}\left(x_{j}\left(t\right)\right)-f_{i}\left(x^{*}\right)\right]\nn\\
&\leq&\frac{1}{T}\sum_{t=1}^{T}\frac{1}{m}\sum_{i=1}^{m}\left[f_{i}\left(x_{j}\left(t\right)\right)-f_{i}\left(y\left(t\right)\right)\right] \nn\\
&&+\frac{1}{T}\sum_{t=1}^{T}\frac{1}{m}\sum_{i=1}^{m}\left\langle g_{i}\left(t\right),y\left(t\right)-x^{*}\right\rangle.
\label{theorem-basic-3}
\EEAS
Using the result of Lemma \ref{lemma-gxxstar} and $\left\Vert g_{i}\left(t\right)\right\Vert_*\leq L$,
\BEAS
&&\frac{1}{T}\sum_{t=1}^{T}\frac{1}{m}\sum_{i=1}^{m}\left\langle g_{i}\left(t\right),y\left(t\right)-x^{*}\right\rangle\nn\\ 
&\leq&\frac{L^{2}}{2T}\sum_{t=1}^{T}\alpha\left(t\right)+\frac{1}{T\alpha\left(T\right)}\psi\left(x^{*}\right).
\label{theorem-basic-4}
\EEAS
Applying the assumption that for $i\in\left\{1,\ldots,m\right\}$, each function $f_{i}$ is Lipstchitz continuous,
\BEAS
f_{i}\left(y\left(t\right)\right)-f_{i}\left(x_{j}\left(t\right)\right)
&\leq& L\left\Vert x_{j}\left(t\right)-y\left(t\right)\right\Vert
\label{theorem-basic-5}
\EEAS
Then using Lemma \ref{lemma-uv}, for $t\geq1$ it achieves that
\BEAS
\left\Vert x_{i}\left(t\right)-y\left(t\right)\right\Vert
\leq \alpha\left(t\right)\left\Vert \frac{z_{i}\left(t\right)}{w_{i}\left(t\right)}-\bar{z}\left(t\right)\right\Vert _{*}.
\label{theorem-basic-6}
\EEAS
We finally take formulas (\ref{theorem-basic-2})-(\ref{theorem-basic-6}) into the inequality (\ref{theorem-basic-1}), then the result of the inequality (\ref{theorem-basic-0}) is carried out. The proof is complete.
\hfill$\blacksquare$

Under the foundation of Theorem \ref{theorem-basic}, we then begin to analyze establishment of Theorem \ref{theorem-alpha}.

{\it Proof of Theorem \ref{theorem-alpha}:}

Considering the term $\left\Vert \frac{z_{i}\left(t\right)}{w_{i}\left(t\right)}-\bar{z}\left(t\right)\right\Vert_* $, with the algorithm (\ref{algorithm-2}) and the definition of  $\bar{z}\left(t\right)$, for $t\geq1$,
\BEAS
\frac{z_{i}\left(t\right)}{w_{i}\left(t\right)}-\bar{z}\left(t\right)
&=&\frac{\sum_{s=0}^{t-1}\sum_{j=1}^{m+\tau}\left[Q\left(t-1:s+1\right)\right]_{ij}g_{j}\left(s\right)}{\sum_{r=1}^{m+\tau}\left[Q\left(t-1:0\right)\right]_{ir}w_j\left(0\right)}\nn\\
&&-\sum_{s=0}^{t-1}\frac{1}{m}\sum_{j=1}^{m+\tau}g_{j}\left(s\right).
\label{theorem-alpha_1}
\EEAS
For $j\in\left\{m+1,\ldots,m+\tau\right\}$ and all the $t$, $g_j\left(t\right)=0$ and $w_j\left(t\right)=0$; and for $j\in\left\{1,\ldots,m\right\}$ and all the $t$, $w_j\left(t\right)=1$. Then the following equation can be presented by subtracting and adding the term $\sum_{s=0}^{t-1}\sum_{j=1}^{m}\phi_{i}\left(t-1\right)g_{j}\left(s\right)$ on the numerator of the right hand side,
\BEAS
\frac{z_{i}\left(t\right)}{w_{i}\left(t\right)}-\bar{z}\left(t\right)
=\frac{\Upsilon_{1}+\Upsilon_{2}}{\sum_{r=1}^{m}\left[Q\left(t-1:0\right)\right]_{ir}},
\label{theorem-alpha_2}
\EEAS
where
\BEASN
&&\Upsilon_{1}=\sum_{s=0}^{t-1}\sum_{j=1}^{m}\left(\left[Q\left(t-1:s+1\right)\right]_{ij}-\phi_{i}\left(t-1\right)\right)g_{j}\left(s\right)\nn\\ &&\Upsilon_{2}=\sum_{s=0}^{t-1}\sum_{j=1}^{m}\left(\phi_{i}\left(t-1\right)-\frac{1}{m}\sum_{r=1}^{m}\left[Q\left(t-1:0\right)\right]_{ir}\right)g_{j}\left(s\right).
\EEASN
Using Corollary \ref{Collary-bound} and the fact that $\left\Vert g_{j}\left(s\right)\right\Vert _{*}\leq L$, for the term $\Upsilon_{1}$ we have,
\BEAS
\left\Vert \Upsilon_{1}\right\Vert _{*}\leq\sum_{s=0}^{t-1}\sum_{j=1}^{m}C\lambda^{t-s-2}\left\Vert g_{j}\left(s\right)\right\Vert _{*}\leq\frac{mCL}{\left(1-\lambda\right)\lambda}.
\label{theorem-alpha_3}
\EEAS
For the other term $\Upsilon_{2}$, it is obtained that
\BEAS
\left\Vert \Upsilon_{2}\right\Vert _{*}
&\leq&\sum_{s=0}^{t-1}\sum_{j=1}^{m}\frac{1}{m}\left(mC\lambda^{t-1}\right)L\nn\\
&\leq&mCt^{*}\lambda^{t^{*}-1}L,
\label{theorem-alpha_4}
\EEAS
where $t^{*}$ is the variable to get the maximum value of the function $t\lambda^{t-1}$ (\ie, $t^{*}=\arg\max\left(t\lambda^{t-1}\right)$). Using Corollary \ref{Collary-bound}, it is known that $\sum_{j=1}^{m}\left[Q\left(t-1:0\right)\right]_{ij} \geq \delta$, and with formulas (\ref{theorem-alpha_2})-(\ref{theorem-alpha_4}), we can give a bound to the term $\left\Vert \frac{z_{i}\left(t\right)}{w_{i}\left(t\right)}-\bar{z}\left(t\right)\right\Vert _{*}$,
\BEAS
\left\Vert \frac{z_{i}\left(t\right)}{w_{i}\left(t\right)}-\bar{z}\left(t\right)\right\Vert _{*}&\leq&\frac{\left\Vert \Upsilon_{1}\right\Vert _{*}+\left\Vert \Upsilon_{2}\right\Vert _{*}}{\delta}\nn\\
&\leq&\frac{\frac{mC}{\left(1-\lambda\right)\lambda}L+mCt^{*}\lambda^{t^{*}-1}L}{\delta}.
\label{theorem-alpha_5}
\EEAS
Finally, with the bound of the term $\left\Vert \frac{z_{i}\left(t\right)}{w_{i}\left(t\right)}-\bar{z}\left(t\right)\right\Vert _{*}$ in the inequality (\ref{theorem-alpha_5}), applying Theorem \ref{theorem-basic}, we obtain the inequality (\ref{theorem-alpha-0}). We have completed the proof.
\hfill$\blacksquare$

\section{Simulation examples}
In this section, we would like to investigate the performance of PS-DDA algorithm with communication delays in solving a distributed quadratic programming problem and a distributed estimation problem in the sensor network.

\subsection{Distributed quadratic programming}\label{section_s1}
We consider a network with $m$ nodes and the each node has its own quadratic function. the distributed quadratic programming problem can be expressed as
\BEASN
\min\frac{1}{m}\sum_{i=1}^{m}\left(x-U_{i}\right)^{T}\left(x-U_{i}\right),\quad\left\Vert x\right\Vert _{1}\leq3,
\EEASN

where $x\in\mathbb{R}^{2}$, and $U_{i}\in\mathbb{R}^{2}$ is a fixed vector known by node $i$. In this example, the network consists of $m=8$ nodes, and each node in the system possesses its own quadratic function, being able to communicate with its neighbors through a union of directed graphs. Four directed graphs $\mathcal{G}_1$, $\mathcal{G}_2$, $\mathcal{G}_3$ and $\mathcal{G}_4$ can be represented by column stochastic weight matrices $Q_1$, $Q_2$, $Q_3$ and $Q_4$ respectively. The weight matrices are given in the following:
\BEASN
&&Q_{1}=\begin{bmatrix}
\begin{smallmatrix}
0.5 & 0 & 0 & 0 & 0   & 0 & 0 & 0\\
0   & 1 & 0 & 0 & 0   & 0 & 0 & 0\\
0.5 & 0 & 1 & 0 & 0   & 0 & 0 & 0\\
0   & 0 & 0 & 1 & 0   & 0 & 0 & 0\\
0   & 0 & 0 & 0 & 0.5 & 0 & 0 & 0\\
0   & 0 & 0 & 0 & 0   & 1 & 0 & 0\\
0   & 0 & 0 & 0 & 0   & 0 & 1 & 0\\
0   & 0 & 0 & 0 & 0.5 & 0 & 0 & 1
\end{smallmatrix}
\end{bmatrix},\
Q_{2}=\begin{bmatrix}
\begin{smallmatrix}
1 & 0   & 0 & 0   & 0 & 0.5 & 0 & 0\\
0 & 0.5 & 0 & 0   & 0 & 0   & 0 & 0\\
0 & 0   & 1 & 0   & 0 & 0   & 0 & 0\\
0 & 0   & 0 & 0.5 & 0 & 0   & 0 & 0\\
0 & 0.5 & 0 & 0   & 1 & 0   & 0 & 0\\
0 & 0   & 0 & 0   & 0 & 0.5 & 0 & 0\\
0 & 0   & 0 & 0.5 & 0 & 0   & 1 & 0\\
0 & 0   & 0 & 0   & 0 & 0   & 0 & 1
\end{smallmatrix}
\end{bmatrix},\\
&&Q_{3}=\begin{bmatrix}
\begin{smallmatrix}
1 & 0   & 0 & 0 & 0 & 0 & 0 & 0\\
0 & 0.5 & 0 & 0 & 0 & 0 & 0 & 0\\
0 & 0   & 1 & 0 & 0 & 0 & 0 & 0\\
0 & 0.5 & 0 & 1 & 0 & 0 & 0 & 0\\
0 & 0   & 0 & 0 & 1 & 0 & 0 & 0\\
0 & 0   & 0 & 0 & 0 & 1 & 0 & 0.5\\
0 & 0   & 0 & 0 & 0 & 0 & 1 & 0\\
0 & 0   & 0 & 0 & 0 & 0 & 0 & 0.5
\end{smallmatrix}
\end{bmatrix},\
Q_{4}=\begin{bmatrix}
\begin{smallmatrix}
1 & 0 & 0   & 0 & 0.5 & 0 & 0   & 0\\
0 & 1 & 0.5 & 0 & 0   & 0 & 0   & 0\\
0 & 0 & 0.5 & 0 & 0   & 0 & 0   & 0\\
0 & 0 & 0   & 1 & 0   & 0 & 0   & 0\\
0 & 0 & 0   & 0 & 0.5 & 0 & 0   & 0\\
0 & 0 & 0   & 0 & 0   & 1 & 0.5 & 0\\
0 & 0 & 0   & 0 & 0   & 0 & 0.5 & 0\\
0 & 0 & 0   & 0 & 0   & 0 & 0   & 1
\end{smallmatrix}
\end{bmatrix}.
\EEASN

We first plot the simulation result of function errors versus iteration span $T$ of the 8 nodes with a communication delay $\tau_{edge}=4$ of each edge in the network in Figure 2. It can be easily found the PS-DDA algorithm with communication delays could converge to the optimum value with time.

\begin{figure}[H]

\begin{center}
\rotatebox{360}{\scalebox{0.4}[0.4]{\includegraphics{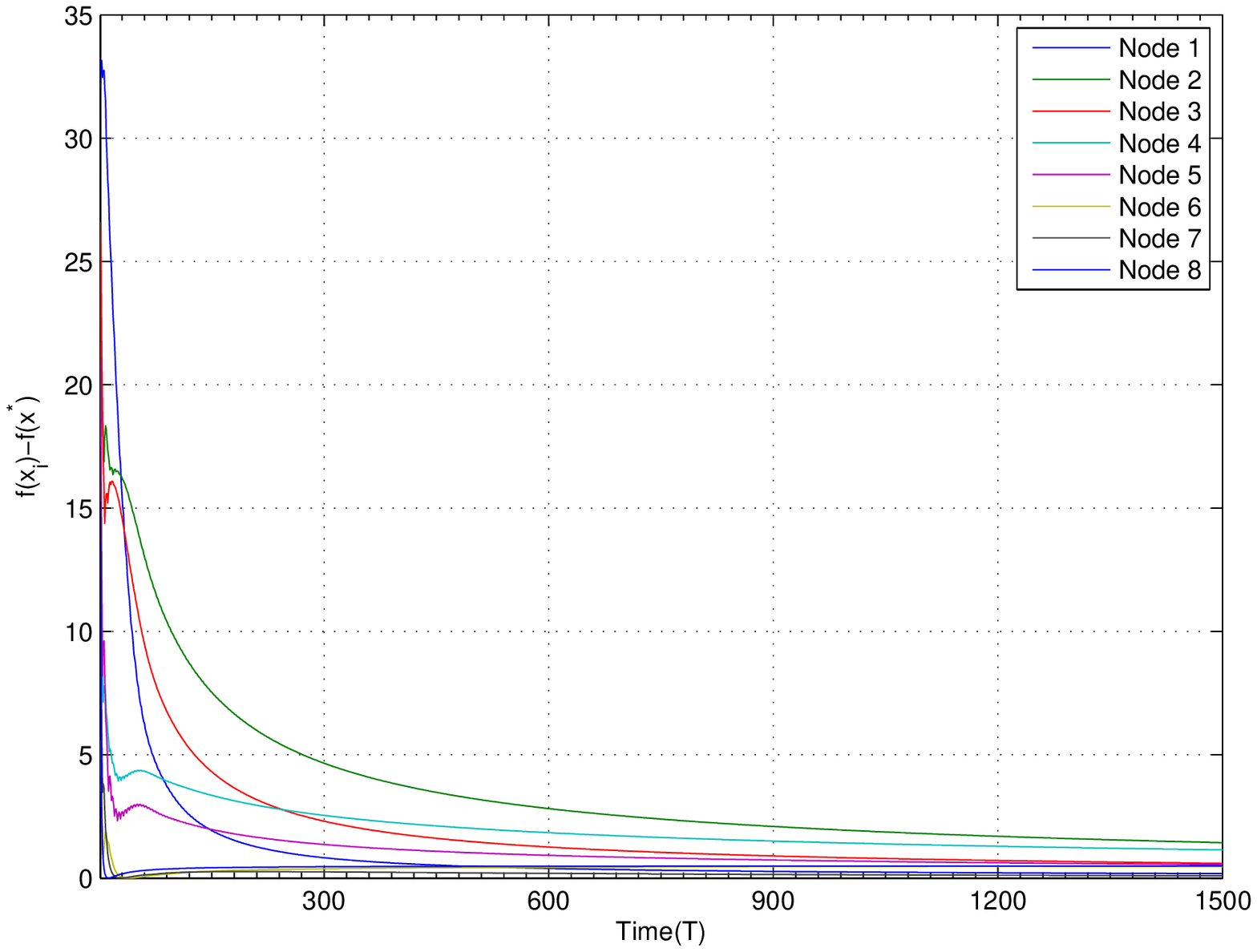}}}
\\[0pt]
{\scriptsize {\vc{Figure 2. The function errors versus iteration span $T$ of the 8 nodes for the distributed quadratic programming problem.} }}
\end{center}
\end{figure}

We then want to make some comparisons with the results of the PS-DDA with communication delays, the distributed dual averaging (DDA) algorithm with communication delays possessing doubly stochastic weight matrix in \cite{KI2012}, and the DDA algorithm with communication delays possessing column stochastic weight matrix. We propose the maximum function errors versus iteration span $T$ of these algorithms in Figure 3. It can be noticed that PS-DDA algorithm with communication delays performs no worse than that of DDA algorithm. Moreover, the PS-DDA algorithm with communication delays does not require doubly stochastic weight matrix and the knowledge of graph sequence and size of graphs.

\begin{figure}[H]

\begin{center}
\rotatebox{360}{\scalebox{0.4}[0.4]{\includegraphics{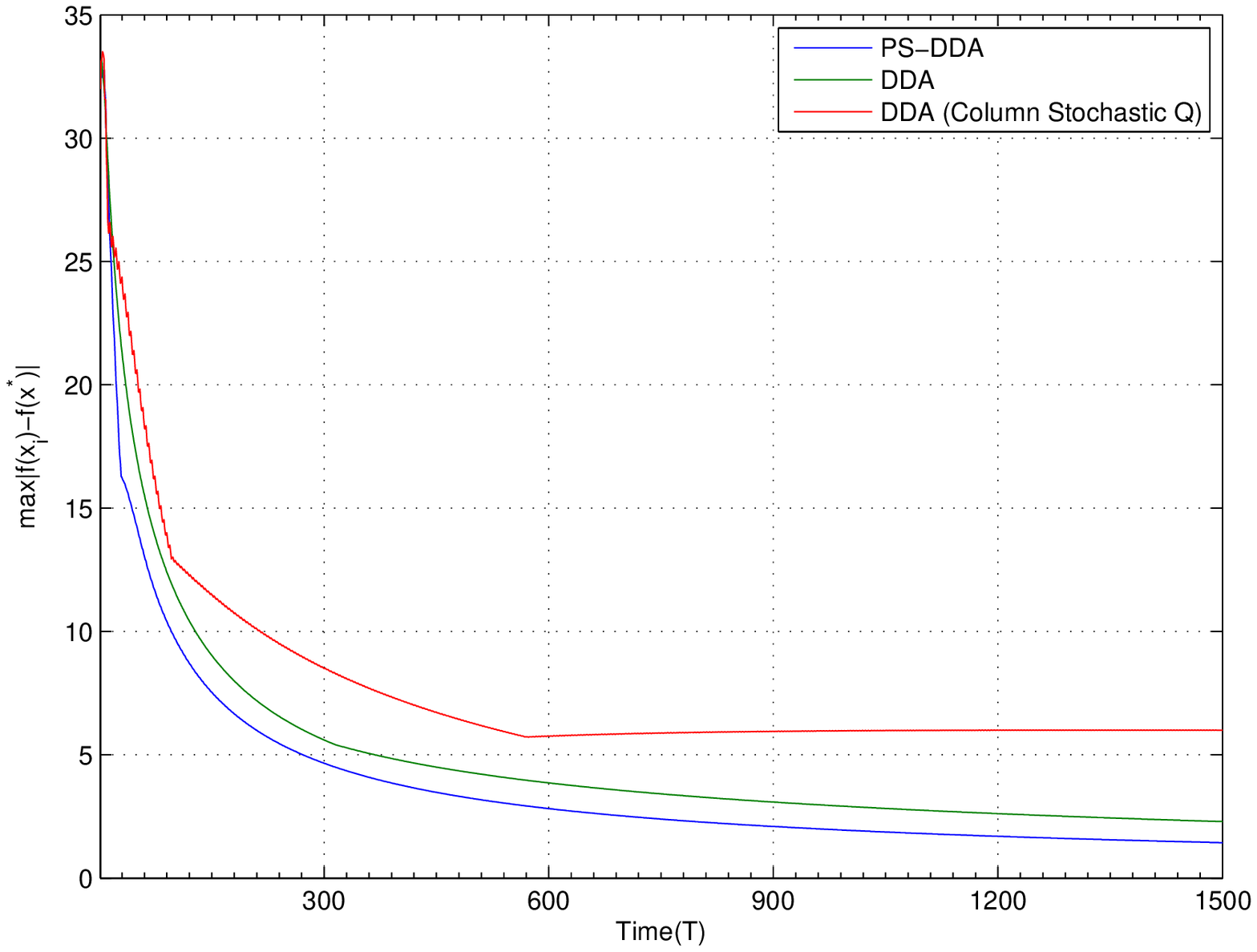}}}
\\[0pt]
{\scriptsize{\vc{Figure 3. The maximum function errors versus iteration span $T$ of 3 algorithms for the distributed quadratic programming problem.} }}
\end{center}
\end{figure}

Finally, the simulation results of the maximum function errors versus iteration span $T$ possessing disparate communication delays of each edge are presented in Figure 4, in which we set the communication delay of each edge as $\tau_{edge}=0$, $\tau_{edge}=4$ and $\tau_{edge}=8$. Obviously, the algorithm with longer the communication delay is, the worse the convergence result will be. The simulation results presented in Figure 4 satisfy Theorem \ref{theorem-alpha} in this paper.

\begin{figure}[H]

\begin{center}
\rotatebox{360}{\scalebox{0.4}[0.4]{\includegraphics{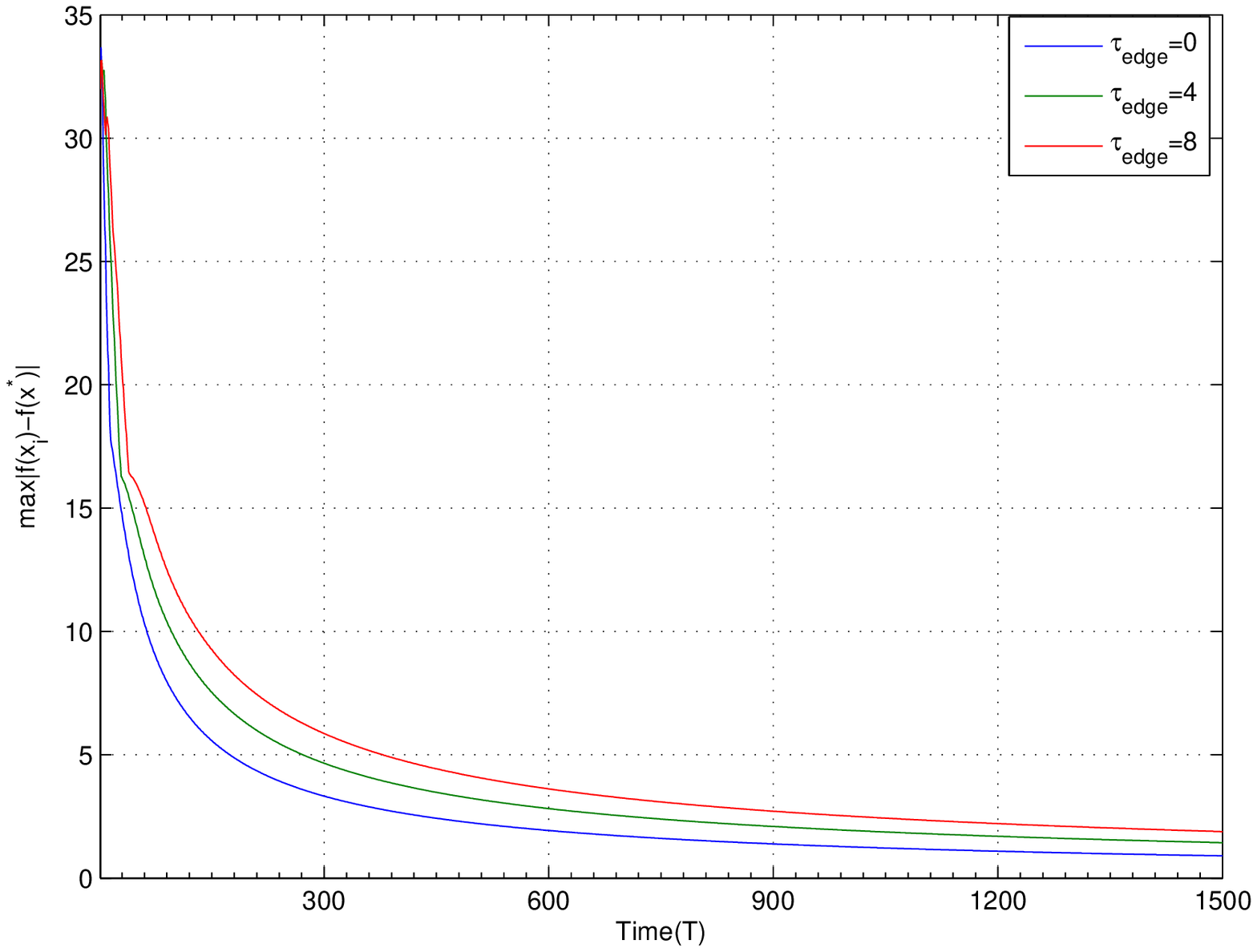}}}
\\[0pt]
{\scriptsize{\vc{Figure 4. The maximum function errors versus iteration span $T$ with disparate communication delays of each edge as $\tau_{edge}=0$, $\tau_{edge}=4$ and $\tau_{edge}=8$ for the distributed quadratic programming problem.} }}
\end{center}
\end{figure}

\subsection{Distributed estimation problem in the sensor network}

In this subsection, the distributed estimation problem in the sensor network is considered (Motivated by \cite{Akbari2017}). Each sensor $i$ in the network system observe a vector $x\in\mathbb{R}^{d}$ through an observation vector $R_{i}\left(x\right)\in\mathbb{R}^{d_i}$. Every sensor is modeled as $S_{i}\left(x\right)=K_{i}x$, where $K_{i}\in\mathbb{R}^{d_i\times d}$. The purpose of the network system is to find the argument $\hat{x}$ to minimize the following function.
\BEASN
\min\frac{1}{m}\sum_{i=1}^{m}\frac{1}{2}\left\Vert R_{i}\left(x\right)-K_{i}\hat{x}\right\Vert ^{2},\quad\left\Vert \hat{x}\right\Vert _{1}\leq h,
\EEASN
where $h$ is positive constant. In this example we utilize the same sequence of network graphs and weight matrices $Q_1$, $Q_2$, $Q_3$ and $Q_4$ in Section \ref{section_s1}. The dimension of the observation vector is $d=1$, and $K_{i}=1$ for all $i$. The observation model for sensor $i$ is
\BEASN
R_{i}\left(x\right)=a_{i}x+b_{i},
\EEASN
where $a_i$ and $b_i$ are generated from uniform distributions on $\left[1,2\right]$ and $\left[-\frac{1}{2},\frac{1}{2}\right]$, respectively, and we select $h=0.1$.

We first simulate with the communication delay of $\tau_{edge}=4$ on each edge. Figure 5 shows the function errors versus iteration span $T$ for of 8 sensors in the network for the distributed estimation problem in the sensor network. It can be found PS-DDA algorithm with communication delays converges in this case.

We then study the distributed estimation problem in the sensor network with different selections of communication delays. The problem is simulated with different communication delays on each edge as $\tau_{edge}=0$, $\tau_{edge}=4$, $\tau_{edge}=8$. Figure 6 proposes the maximum function errors versus iteration span $T$ with disparate communication delays of each edge as $\tau_{edge}=0$, $\tau_{edge}=4$ and $\tau_{edge}=8$ for the distributed estimation problem in the sensor network. It also can be found from Figure 6 that the communication delay will affect the convergence of the PS-DDA algorithm.

\begin{figure}[H]

\begin{center}
\rotatebox{360}{\scalebox{0.4}[0.4]{\includegraphics{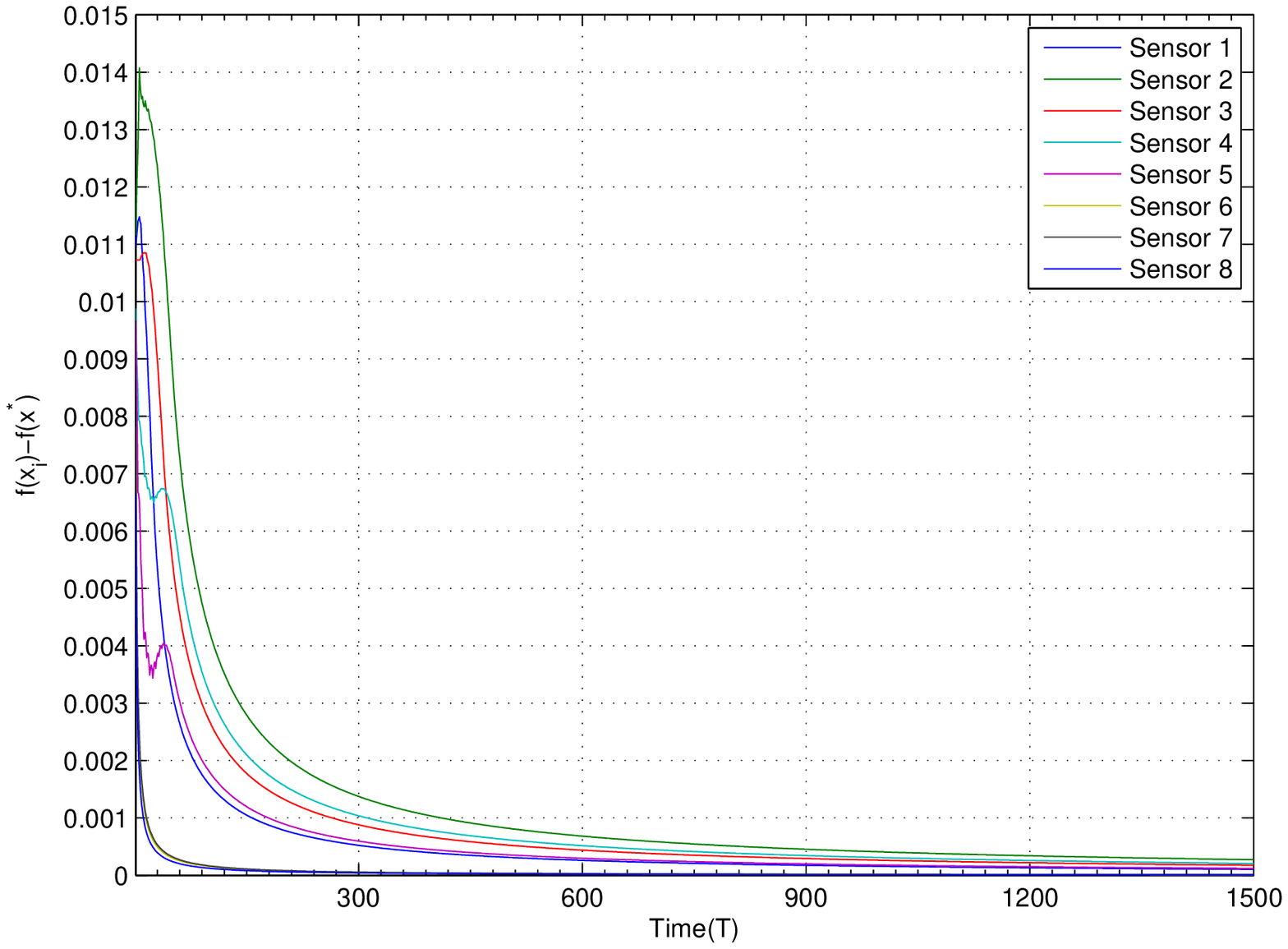}}}
\\[0pt]
{\scriptsize {\vc{Figure 5. The function errors versus iteration span $T$ of 8 sensors in the network for distributed estimation problem in the sensor network.} }}
\end{center}
\end{figure}

\begin{figure}[H]

\begin{center}
\rotatebox{360}{\scalebox{0.4}[0.4]{\includegraphics{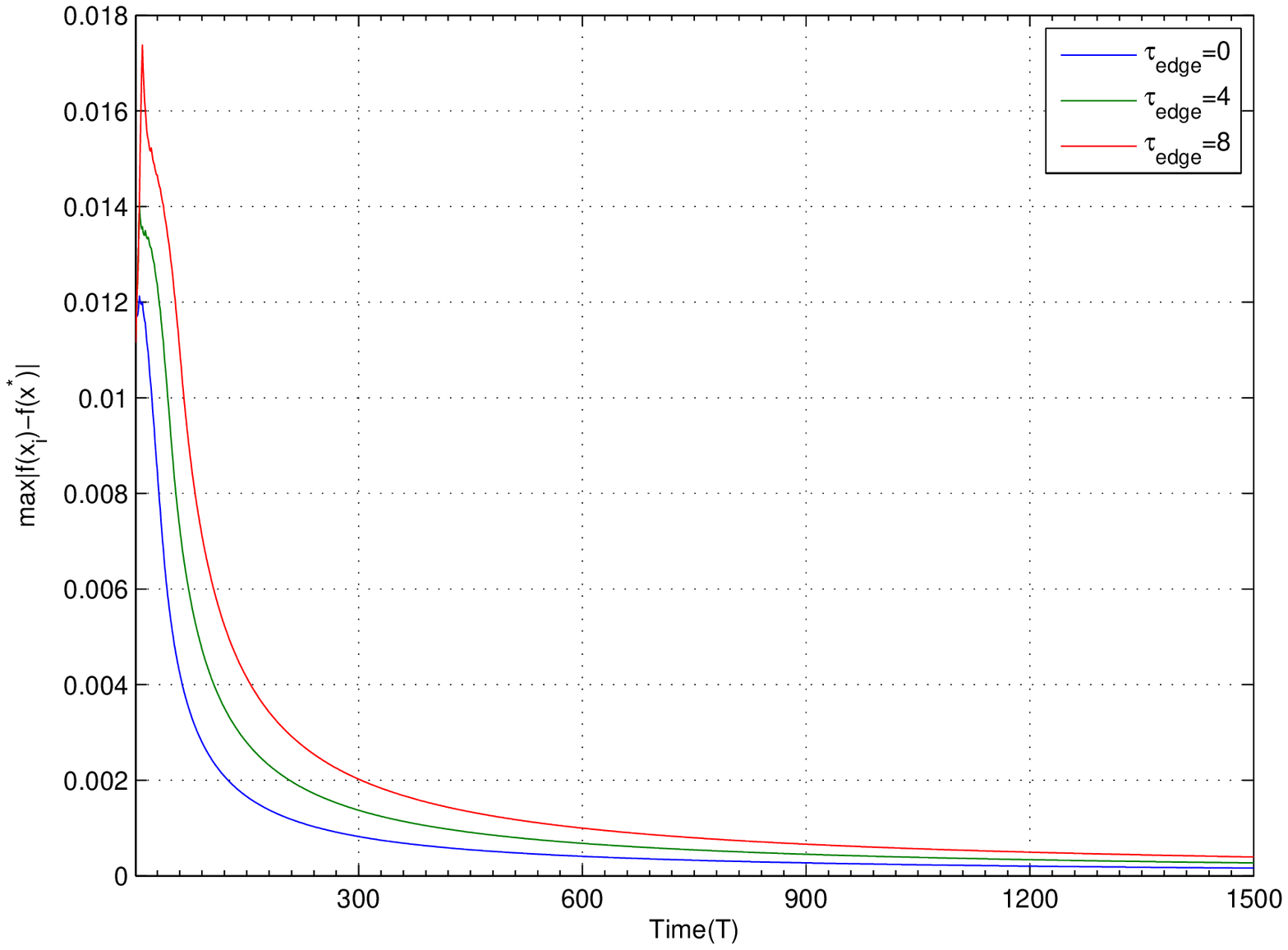}}}
\\[0pt]
{\scriptsize{\vc{Figure 6. The maximum function errors versus iteration span $T$ with disparate communication delays of each edge as $\tau_{edge}=0$, $\tau_{edge}=4$ and $\tau_{edge}=8$ for distributed estimation problem in the sensor network.} }}
\end{center}
\end{figure}

\section{Conclusions}
We analyze PS-DDA algorithm with non-negligible communication delays that the original algorithm ignores. The fixed communication delay model in \cite{KI2011} is applied and PS-DDA algorithm could converge at a rate $\calO\left(T^{-0.5}\right)$ with proper step size, where $T$ is iteration span. Simulations to numerical examples are presented in our paper, and the performance of simulations verify the convergence result of the PS-DDA algorithm with communication delays. Following two aspects can be further explored in the future. On one hand, the convex distributed optimization problem with communication delays could be considered by some other algorithms like the distributed mirror descent algorithm. On the other hand, the communication delays are assumed to be fixed in this technical correspondence, and this is a constraint to application. As a result, considering the PS-DDA algorithm with random communication delays may also be a future research direction.


\end{document}